\documentclass[10pt]{article}
\usepackage[utf8]{inputenc}
\usepackage{mathptmx}
\usepackage[scaled=.90]{helvet}
\usepackage{graphicx}
\usepackage{amssymb}
\usepackage{bbm}
\usepackage[numbers, sort, comma, square]{natbib}
\usepackage[english]{babel}
\usepackage{amsmath}
\usepackage{amsthm}
\usepackage{thmtools, thm-restate}
\usepackage{hyperref}
\usepackage{cleveref}
\usepackage{comment}
\usepackage{geometry}
\usepackage{tikz}
\usepackage{tkz-tab}
\usetikzlibrary{shapes}
\usepackage{floatrow}
\usepackage{mathtools}
\usepackage{wrapfig}
\usepackage{gensymb}
\usetikzlibrary{patterns}
\usetikzlibrary{calc,backgrounds}
\usepackage{stmaryrd}
\usepackage{caption}
\usepackage{subcaption}
\usepackage{fancyhdr}
\usepackage{float}
\usepackage{rotating}
\usepackage{calc}
\usepackage{lipsum}
\usepackage{accents}
\usepackage{filecontents}
\usepackage{awesomebox}
\usepackage{alertmessage}
\usepackage{titlesec}
\usepackage[toc,page]{appendix}
\usetikzlibrary{arrows.meta}
\tikzset{>={Latex[width=0.5mm,length=0.5mm]}}
\usepackage{pgfplots}
\pgfplotsset{compat=1.15}
\usepackage{mathrsfs}
\usepackage{eso-pic}
\usetikzlibrary{arrows}
\usepackage{tikz-cd}
\usepackage{bm}
\usepackage{bbold}
\usepackage{indentfirst}

\makeatletter
\newcommand{\mscmotscles}[3]{%
  \begingroup
  \renewcommand{\thefootnote}{}
  \renewcommand{\@makefntext}[1]{\noindent ##1}
  \footnotetext{%
    \textit{2020 Mathematics Subject Classification :} #1.\\
    \textit{Keywords :} #2.\\
    #3}%
  \endgroup
}
\makeatother

\newcommand{\adresseauteur}[6]{%
  \par\addvspace{2em}%
  \noindent
  \begin{minipage}{\linewidth}
    #1\\                      
    \textit{#2}\\             
    \textit{#3}\\             
    \textit{#4}\\             
    \textit{#5}\\             
    \href{mailto:#6}{\nolinkurl{#6}}
  \end{minipage}%
  \par
}

\titleformat{\section}[block]{\Large\bfseries}{\thesection.}{1em}{}
\titleformat{\subsection}[block]{\hspace{2em}\large\bfseries}{\thesubsection}{1em}{}
\titleformat{\subsubsection}[block]{\hspace{3em}\large\bfseries}{\thesubsubsection}{1em}{}

\makeatletter
\newcommand\subsubsubsection{\@startsection{paragraph}{4}{\z@}{-2.5ex\@plus -1ex \@minus -.25ex}{1.25ex \@plus .25ex}{\hspace{3em}\normalfont\large\bfseries}}
\newcommand\subsubsubsubsection{\@startsection{subparagraph}{5}{\z@}{-2.5ex\@plus -1ex \@minus -.25ex}{1.25ex \@plus .25ex}{\hspace{4em}\normalfont\large\bfseries}}
\makeatother

\renewcommand*{\thefootnote}{\fnsymbol{footnote}}

\theoremstyle{plain}

\newtheorem{theorem}{Theorem}
\newtheorem*{theorem*}{Theorem}
\newtheorem{proposition}[theorem]{Proposition}
\newtheorem{corollary}[theorem]{Corollary}
\newtheorem*{corollary*}{Corollary}
\newtheorem{lemma}[theorem]{Lemma}
\newtheorem{axiom}[theorem]{Axiom}
\newtheorem{conjecture}{Conjecture}
\newtheorem*{conjecture*}{Conjecture}

\theoremstyle{remark}
\newtheorem{remark}[theorem]{Remark}
\newtheorem{example}[theorem]{Example}

\theoremstyle{definition}
\newtheorem{definition}[theorem]{Definition}

\newenvironment{lem}{\begin{lemma}~\newline}{\end{lemma}}

\makeatother

\numberwithin{equation}{section}
\numberwithin{theorem}{section}
\numberwithin{definition}{section}
\numberwithin{proposition}{section}
\numberwithin{lemma}{section}
\numberwithin{remark}{section}
\numberwithin{axiom}{section}
\numberwithin{corollary}{section}

\let\oldsection\section
\renewcommand{\section}{
  \renewcommand{\theequation}{\thesection.\arabic{equation}}
  \renewcommand{\thetheorem}{\thesection.\arabic{theorem}}
  \renewcommand{\theproposition}{\thesection.\arabic{proposition}}
  \renewcommand{\thedefinition}{\thesection.\arabic{definition}}
  \renewcommand{\theremark}{\thesection.\arabic{remark}}
  \renewcommand{\theaxiom}{\thesection.\arabic{axiom}}
  \renewcommand{\thelemma}{\thesection.\arabic{lemma}}
  \renewcommand{\thecorollary}{\thesection.\arabic{corollary}}
  \oldsection}

\newcommand{\N}{\mathbb{N}}
\newcommand{\Ns}{\mathbb{Z}_{> 0}}

\newcommand{\A}{\mathcal{A}}
\newcommand{\Kappa}{\mathcal{K}}

\renewcommand{\S}{\mathfrak{S}}

\newcommand{\cf}{\mathcal{F}}
\newcommand{\cn}{\mathfrak{N}}
\renewcommand{\P}{\mathcal{P}}
\newcommand{\I}{\mathfrak{I}}
\newcommand{\M}{\mathfrak{M}}
\newcommand{\m}{\mathfrak{m}}

\newcommand{\ncurly}{\mathcal{N}}

\newcommand{\frakN}{\mathfrak{N}}
\newcommand{\kfrak}{\mathfrak{k}}
\newcommand{\Kfrak}{\mathfrak{K}}

\renewcommand{\a}{\alpha}
\renewcommand{\b}{\beta}

\newcommand{\bigoh}[1]{\mathcal{O}\left(#1\right)}

\newcommand{\et}{\text{ and }}
\newcommand{\ou}{\text{ or }}
\newcommand{\si}{\text{ if }}

\renewcommand{\th}{\text{th}}

\newcommand{\half}{\frac{1}{2}}
\newcommand{\minus}{{-1}}

\newcommand{\vertiii}[1]{{\left\vert\kern-0.25ex\left\vert\kern-0.25ex\left\vert #1 
    \right\vert\kern-0.25ex\right\vert\kern-0.25ex\right\vert}}

\newcommand{\floor}[1]{\left\lfloor #1 \right\rfloor}

\newcommand{\set}[1]{\left\{ #1 \right\}}

\usepackage{scalerel}
\usepackage{stackengine}
\usepackage{wasysym}

\makeatletter
\def\iddots{\mathinner{\mkern1mu\raise\p@
\vbox{\kern7\p@\hbox{.}}\mkern2mu
\raise4\p@\hbox{.}\mkern2mu\raise7\p@\hbox{.}\mkern1mu}}
\makeatother

\newcounter{Centers}
\newcommand{\leftcentersright}[4][2]{%
	\ifthenelse{\equal{\value{Centers}}{0}}{%
		\leftcentersrightonpage[#1]{#2}{#3}{#4}%
	}{
		\leftcentersrightonitems[#1]{#2}{#3}{#4}%
	}
}

\newlength{\leftlength}
\newlength{\rightlength}
\newcommand{\leftcentersrightonitems}[4][2]{%
	\settowidth{\leftlength}{#2}%
	\settowidth{\rightlength}{#4}%
	\par\calculvskip{#1}%
	\noindent#2\hskip-\leftlength%
	\hfill{#3}\hfill%
	\mbox{}\hskip-\rightlength#4%
	\par\vskip3mm%
}

\newlength{\calculskip}
\newcommand{\calculvskip}[1]{%
	\ifthenelse{#1 = 0}{}{%
		\ifthenelse{#1 = 1}{\setlength{\calculskip}{\smallskipamount}}{}%
		\ifthenelse{#1 = 2}{\setlength{\calculskip}{\medskipamount}}{}%
		\ifthenelse{#1 = 3}{\setlength{\calculskip}{\bigskipamount}}{}%
		\ifthenelse{#1 = 4}{\setlength{\calculskip}{1cm}}{}%
		\ifvmode\vpush{\calculskip}\fi%
	}%
} 

\newcommand{\leftcentersrightonpage}[4][2]{%
	\settowidth{\leftlength}{#2}%
	\settowidth{\rightlength}{#4}%
	\par\calculvskip{#1}
	\noindent#2\hskip-\leftlength%
	\hfill{#3}\hfill
	\mbox{}\hskip-\rightlength#4%
	\par\vskip3mm%
}

\newcommand{\leftcenters}[3][2]{\leftcentersright[#1]{#2}{#3}{}}

\title{\textsc{Representation by sums of even unlike powers}}
\author{\textsc{LO\"{I}C MERMET}}
\date{}

\begin{document}
\maketitle
\mscmotscles
  {Primary 11D85; Secondary 11P55}
  {Fre\u{ı}man's theorem, Hardy-Littlewood method}
  {The author thanks the University of Chalmers for hosting him during periods during which progress was made on this paper, while being supported by Julia Brandes grant's 2022-03717 from Vetenskapsr{\aa}det (Swedish Science Foundation).}

\begin{abstract}
    Fre\u{ı}man's theorem shows that there exists a $s \in \N$ such that all large integers are represented in the form $n = x_1^2+x_2^4+\cdots+x_{s}^{2s}$, with $x_i \in \N$ since the Harmonic series diverges. We show that $79$ powers suffice, improving a previous result of Kuan, Lesesvre and Xiao \cite{KLX2022}. We also adapt this method to provide a qualitative improvement of the result of Liu and Zhao \cite{LiuZhao2021_RepresentationUnlikePowers} by showing that
    all large integers $n$ can be represented in the form
    \vspace*{-10pt}
    \begin{equation*}
        x_1^2 + x_2^3 + \cdots + x_{11}^{12} + x_{12}^{13} + x_{13}^{17} = n.
    \end{equation*}
\end{abstract}
\section{Introduction}

Recent advances in smooth number technology have made possible numerous improvements in results regarding Waring's problem (see \cite{BrudernWooley2023Waringlarge}).
A closely related problem is to represent an integer as the sum of unlike powers. This problem, now known as Fre\u{ı}man-Scourfield's theorem, is the following. Given natural numbers $k_i, \ (i \ge 1)$ satisfying $2 \le k_1 \le k_2 \le \dots$ and $j \in \N$, is there a natural number $s = s(j)$ such that all large integers $n$ can be represented in the form
    \begin{equation}
        \label{Freiman}
        x_1^{k_j}+x_2^{k_{j+1}}+\cdots + x_s^{k_{j+s-1}} = n,
    \end{equation}
    with $x_i$ being some non-negative integers? 
In 1949, Fre\u{ı}man announced in \cite{Freiman1949aClaim} that such a result would hold as long as $\sum \frac{1}{k_j}$ diverges, and he provided a sketch of a proof in the same year \cite{Freiman1949bVariantWaring}. A decade later, Scourfield \cite{Scourfield1963GeneralizationWaring} filled in the gaps present in Fre\u{ı}man's original proof. Since then, various authors have made improvements to related problems and have aimed to give an unconditional bound on $s$ (see \cite{BrudernWooley2023Freiman,Kufner2025}).
For instance, when $(k_j)_j$ lies in an arithmetic progression, various works cover this particular case. Indeed, let $k_j = k_0+r(j-1)$ for some $k_0 \in \Ns, \ r\in \Ns$ and $j \in \N$. Then, following Ford \cite{Ford1996UnlikePowersII}, we denote by $R(k_0,r)$ the least number of variables $s$ such that \eqref{Freiman} holds for sufficiently large integers with the set of exponents $k_j$. The most notable values of $k_0 \et r$ are obviously $(k_0,r) = (k,0)$ as one recovers the original Waring's problem. The second interesting non-trivial couple is $(k_0,r) = (2,1)$. Decades of work, initiated by Roth (e.g. \cite{Roth1951AdditiveNumberTheory,Thanigasalam1983CertainAdditiveRepresentations,Ford1995UnlikePowersI,Ford1996UnlikePowersII,LiuZhao2021_RepresentationUnlikePowers}) led Liu and Zhao in 2021 to finally achieve that
\begin{equation}
    \label{UnlikePowers}
    R(2,1) \le 13.
\end{equation}
Ford, in \cite{Ford1996UnlikePowersII}, has also proved that $R(3,1) \le 72$, and in their breathtaking paper, Br\"{u}dern and Wooley have reached $R(k,1) \le 100(k+1)^2$ (see \cite[Corollary 1.4]{BrudernWooley2023Freiman}) and 

$$R(k,r) \le r^\minus25^r(r+1)^{r+1}(k+1)^{r+1},$$ 

\noindent for any integers $k \et r$, and even better $R(k,r) \le (6k+6)^{2r}$ whenever $r \ge k$ (see \cite[Theorem 1.3]{BrudernWooley2023Freiman}).
Recent work carried out by Kuan, Lesesvre and Xiao \cite{KLX2022} claims that, for even exponents, one has $R(2,2) \le 133$. Using more recent technology, we made a significant improvement in that direction to achieve the following.
\begin{theorem}
    \label{ThEvenUnlikePowers}
    All large integers $n$ can be represented in the form
    \begin{equation}
        \label{EqEvenPowers}
        x_1^2 + x_2 ^4 + \cdots + x_{79}^{158} = n
    \end{equation}
    with $x_i \in \N$. In other words, one has
    \begin{equation}
    R(2,2) \le 79.
    \end{equation}
\end{theorem}

Finally, a related problem has been partially solved by Br\"{u}dern in \cite{Brudern2019}, in which he proved that for $s=5$, almost every large integer can be written as the sum of a square, a biquadrate, a sixth power, an eighth power and a tenth power. This result uses very different tools, but also shows that there is still room to improve the previously stated theorem. 

In that direction, we suggest a quantitative improvement of the result of Liu and Zhao about $R(2,1)$ following the approach presented in this paper. Indeed, we ask for the largest integer $s \ge 2$ such that all sufficiently large integers $n$ can be represented as 

\begin{equation}
\label{FreimanSSup}
x_1^2 + x_2^3 + \cdots + x_{11}^{12} + x_{12}^{13} + x_{13}^{s} = n.
\end{equation}

We know that $s = 14$ works, as proved by Liu and Zhao \cite{LiuZhao2021_RepresentationUnlikePowers}. Also, as the density of powers $k$ decreases when $k$ increases, the heuristic emphasises the fact that obtaining $s$ as large as possible would provide a solution to reach $R(2,1) \le 13$. We claim that $s=17$ is attainable, as stated in the following theorem.

\begin{restatable}[]{thm}{ThUnlikePowersSupresta}
\label{ThUnlikePowersSup}
    All large integers $n$ can be represented in the form
    \begin{equation}
        x_1^2 + x_2^3 + \cdots + x_{11}^{12} + x_{12}^{13} + x_{13}^{17} = n
    \end{equation}
    with $x_i \in \N$.
\end{restatable}

We briefly remark that the circle method is expected to work well as long as the sum of the inverses of the exponents is greater than 2. Thus, this quantitative improvement we just provided can be viewed as
$$\sum\limits_{k=2}^{13}\frac{1}{k} \simeq 2.1801 \ ;  \sum\limits_{k=2}^{13}\frac{1}{k}+\frac{1}{17} \simeq 2.2389 \et \sum\limits_{k=2}^{14}\frac{1}{k} \simeq 2.2515$$ 

Again, we will adapt the method used to claim \autoref{ThEvenUnlikePowers} to provide a new bound for $R(3,1)$ as Ford introduced it in \cite{Ford1996UnlikePowersII}.
\begin{restatable}[]{thm}{ThUnlikePowersThreeresta}
    \label{ThUnlikePowersThree}
     All large integers $n$ can be represented in the form
    \begin{equation}
        \label{EqThreePowers}
        x_1^3 + x_2 ^4 + \cdots + x_{71}^{73} = n,
    \end{equation}
    with $x_i \in \N$. In other words, one has
    \begin{equation}
    R(3,1) \le 71.
    \end{equation}
\end{restatable}

The structure of our paper is as follows. Section \ref{SecMVE} will start by providing some insight into mean value estimates for smooth Weyl exponential sums and Ford's method \cite{Ford1996UnlikePowersII} to achieve the first main lemma of the proof.
Section \ref{Pruning} will be devoted to a pruning of the major arcs using some recent results of Br\"{u}dern and Wooley \cite{BrudernWooley2024SmoothMajor}.
Subsections \ref{SecSingularSeries} and \ref{SecSingularIntegral} are dedicated to the study of the singular series and singular integral that naturally appear on such small major arcs; \autoref{ThEvenUnlikePowers} will follow. Then, Section \ref{UAP} shows how to use this method to enhance the original result of \cite{LiuZhao2021_RepresentationUnlikePowers} for $R(2,1)$ via \autoref{ThUnlikePowersSup}. Finally, Section \ref{StartThree} will set out the key ideas for \autoref{ThUnlikePowersThree}.

\section{Notation}
In this paper, the letter $p$ is reserved to denote a prime number. We write $e(z)$ for $e^{2i \pi z}$, and we make extensive use of Vinogradov's notations $\ll$ and $ \gg$. In the following, the letter $\varepsilon$ will denote a sufficiently small positive real number. Remark that any statement in which $\varepsilon$ appears holds for each fixed $\varepsilon >0$, and thus any constants implied in such a statement also depend on $\varepsilon$.

\section{Some infrastructure}
Writing $H(n)$ for the number of solutions to \eqref{EqEvenPowers}, a classic trick using Fourier's orthogonality leads to
$$
H(n) = \int_0^1 f_2(\a)\cdots f_{158}(\a) e(-n\a)\text{d}\a,
$$
where 
\begin{equation}
\label{Deffk}
    f_k(\alpha) = \sum\limits_{x \le X_k} e(x^k \alpha)
\end{equation}
denotes the Weyl sum for some $k \in \N$ greater than or equal to 2, $\a \in [0,1)$ and 
\begin{equation}
    \label{eq:defXk}
    X_k = \floor{n^\frac{1}{k}}.
\end{equation}
 
\noindent However, claiming an asymptotic for $H(n)$ is far beyond the reach of current technology. We can compare this problem with the study of $\tilde{G}(k)$ in Waring's problem, where the best known values (see \cite[Paragraph 14]{Wooley2018}) are significantly less sharp than the current records for $G(k)$. In broad terms, the higher the number of variables, the easier it is to reach an asymptotic for $H(n)$. Now, since we are only interested in proving that $R(n)$ is non-zero for large enough $n$, reducing the scope of variables $(x_i)_i$ to some subsets of $\N^{79}$ may bypass the technical limitations of available technologies regarding complete Weyl sums. 

As mentioned in the preamble, we may adopt a smooth number framework as in \cite{BrudernWooley2023Waringlarge,BrudernWooley2023Freiman} to pursue this direction. We start with a brief overview of the notation used, following that used by Br\"{u}dern and Wooley. Let $\A(P,R)$ denote the set of natural numbers that does not exceed $P$ and with no prime factors larger than $R$,
$$
\A(P,R) = \set{x \in \N^*, x \le P : p|x \implies p \le R}.
$$
When $R = P^\eta$, the set of $R$-smooth numbers is full size, which means that $|\A(P,P^\eta)| \gg_\eta P$ (e.g. \cite[Lemma 12.1]{vaughan1981hardy}). 
Let $k \ge 2$ be an integer, $\eta >0$ and recalling \eqref{eq:defXk}, we define the smooth Weyl sum
\begin{equation}
\label{Defgk}
g_k(\alpha) = \sum\limits_{x \in \A(X_k,X_k^{\eta})}e(x^k\alpha).    
\end{equation}
Then, we introduce
$$\cf = f_2g_{4}g_{6}\cdots g_{158}.$$
A quick study of the integral
\begin{equation}
    \label{Teta}
    T(n;\eta) = \int_0^1 \cf(\a)e(-n\a)\text{d}\a
\end{equation}
shows that $H(n) \ge T(n;\eta)$. Therefore, we may turn our attention to the study of this new integral formulation. As we shall see, the Hardy-Littlewood method delivers an asymptotic bound for $T(n,\eta)$ and actually yields
$$
H(n) \gg \cf(0)n^\minus,
$$
which suffices to prove \autoref{ThEvenUnlikePowers}.

\subsection{Farey dissection}

To successfully apply the circle method to our problem, it is necessary to introduce some infrastructure regarding the Farey dissection employed. Write $L = \log n$ and let $1 \le Q \le \frac{1}{4}\sqrt{n}$ be a parameter. We define

$$\M_{q,a}(Q) := \left\{\a \in [0,1) : \left|a - \frac{a}{q}\right| \le \frac{Q}{nq}\right\}$$

\noindent with $1 \le a \le q \le Q \et \gcd(a,q) = 1.$ Then we define $\M(Q)$ the set of major arcs as

$$\M(Q) = \bigcup\limits_{1 \le q \le Q}\bigcup\limits_{\substack{1 \le a \le q \\ (a,q)=1}} \M_q(Q).$$

\noindent We then define the set of minor arcs $\m(Q)$ as its complement in $[0;1)$. From now on, whenever $\a \in \M_{q,a}(Q) \subset \M(Q)$ is taken, we will use the notation $\b := \a - \frac{a}{q}$ with associated $a \et q$ for $\a$. Remark in particular that for $Q$ sufficiently small (to be precise $Q \le \frac{1}{4}n^\half$), none of the $\M_{q,a}(Q)$ overlaps.
Another interesting set is given by the intermediate arcs $\frakN(Q)$ defined by
$$\frakN(Q) = \M(Q) \backslash \M(Q/2).$$

Introducing such dyadically truncated sets of major arcs is mostly motivated by the observation that they can be seen as major arcs $\M(Q)$ or minor arcs $\m(Q/2)$. Indeed, for $\a \in \cn(Q)$, $1\le a \le q \le Q \et (a,q) = 1$ associated, either
\vspace{-2pt}
$$\frac{Q}{2} < q \le Q \et \left|\a-\frac{a}{q}\right| \le \frac{Q}{nq},$$ 

\leftcenters[0]{or}{$1 \le q \le Q/2 \et \frac{Q}{2nq} \le \left|\a-\frac{a}{q}\right| \le \frac{Q}{nq} .$}

\noindent In both cases, we have that \begin{equation}
    \label{EstNQ} q\left(1+n\left|\a-\frac{a}{q}\right|\right) \asymp Q.
\end{equation}

\noindent The full power of this intermediate dissection will become apparent in Section \ref{MajArcsasymptotic}. Finally, in keeping with the philosophy of the circle method, we expect to prove
$$\int_{\m(Q)} |\cf(\a)|\text{d}\a = o(\cf(0)n^{-1}),$$
with $Q$ as small as possible to be able to easily handle the remaining integral over the set of major arcs. 
More precisely, we aim to reach a logarithmic-sized major arc so that the study of the main term is routine. The following lemma embodies this idea and will be proved in Section \ref{Pruning}.

\begin{lem}
    \label{LemmaMajArcLog}
    We have
    \vspace{-2pt}
    $$T(n;\eta) = \int_{\M(L^{1/15})} f_2^*(\a)g_4(\a)g_6(\a)\cdots g_{158}(\a)\text{d}\a +\bigoh{\cf(0)n^\minus L^{-\rho}},$$
    
    \noindent for some $\rho >0$ and with $f_2^*$ defined below with \eqref{DefFk*}.
\end{lem}

\subsection{Some estimates for classical Weyl sums}

In light of the historical development towards Waring's problem for $k$ large, we introduce results about exponential Weyl sums. Taking $k \in \N^* $ and set $X_k := \floor{n^\frac{1}{k}}$, we write

$$S_k(q,a) = \sum\limits_{x=1}^q e(ax^k/q).$$

\noindent Therefore, for $(a,q)=1$,  \cite[Lemma 2.3]{vaughan1981hardy} states that
\begin{equation}
\label{EstSk}
S_k(q,a) \ll q^{1-\frac{1}{k}+\varepsilon}.
\end{equation}

\noindent Scholars may argue that a slightly stronger estimate for $S_k(q,a)$ can be found in the literature, but we will not need such an improvement in our presentation. In that way, when $\a \in \M(Q)$, we write $\a$ as

$$\a = \frac{a}{q}+\b \quad \text{with } 1\le a \le q \le Q, (a,q)=1 \et |\b| \le \frac{Q}{qn}.
$$
\noindent Finally, define
\begin{equation}
\label{DefFk*}
f_k^*(\a) = \frac{1}{q}S_k(q,a)v_k(\b) \text{ where } \nu_k(\b) = \frac{1}{k}\sum\limits_{m \le n} m^{\frac{1}{k}-1}e(\b m).
\end{equation}

\noindent When it comes to estimate $\nu_k$, \cite[Lemma 2.8]{vaughan1981hardy} produces
\begin{equation}
\label{Estnuk}
v_k(\b) \ll X_k(1+n|\b|)^{-\frac{1}{k}}.
\end{equation}

\noindent Consequently, combining \eqref{EstSk} and \eqref{Estnuk}, we find

$$
f_k^*(\a) \ll X_k q^{-\frac{1}{k}+\varepsilon}(1+n|\b|)^{-\frac{1}{k}}  \quad \text{for} \ \a =\frac{a}{q} + \b.
$$

Lastly, we introduce the well-known Weyl inequality, as it will be of great help to temper the contribution of~$f_2$. Interestingly, note that this inequality will be used only for $f_2$ in Section \ref{FirstMajArc} and $f_3$ in Section \ref{StartThree}. A far better result exists when $k\ge 6$ (see \cite[Theorem 5]{Bou16x}), but will be of no use here, since we restricted our problem to smooth exponential sums when $k$ turns out to be large.
\begin{lemma}{\textit{(Weyl inequality, \cite[Lemma 2.4]{vaughan1981hardy})}}
\label{WeylEstimate}
~\newline
Suppose that $1 \le q\le n^\half, (a,q)=1, |\a - a/q| < q^{-2}$, and let
$$\phi(x) = \a x^k+\a_1x^{k-1}+\cdots+\a_{k-1}x+\a_k.$$
Then
\begin{equation}
    \label{WeylEstimateEq}
    \sum\limits_{x=1}^X e(\phi(x)) \ll X^{1+\varepsilon}\left(\frac{1}{q}+\frac{1}{X}+\frac{q}{X^k}\right)^
    {2^{1-k}}.
\end{equation}
\end{lemma}

\noindent As we shall see later, whenever $\a \in \m(Q)$ for some $1 \le Q \le n^\half$, this lemma yields
\begin{equation}
\label{EstF2Minor}
\sup\limits_{\a \in \m(Q)} |f_2(\a)| \ll f_2(0)^{1+\varepsilon}Q^{-\half}.  
\end{equation}

\noindent Similarly, for any $\a \in \m(Q)$ where $1 \le Q \le n^\frac{1}{3}$, we get

\begin{equation}
\label{EstF3Minor}
\sup\limits_{\a \in \m(Q)}|f_3(\a)|  \ll f_3(0)^{1+\varepsilon}Q^{-\frac{1}{4}}.
\end{equation}

\subsection{Mean value estimates for smooth Weyl sums}
\label{SecMVE}
The smooth Weyl sums $g_k$ that are the key to our results are defined for any $\a \in [0,1)$ by

$$g_k(\a,P,R) = \sum\limits_{x \in \A(P,R)} e(x^k\a).$$

\noindent Define further

\begin{equation}
\label{DefUks}
    U_{k,s}(P) := \int_0^1 |g_k(\a,P,R)|^{s} \text{d}\a  \ll P^{s-k+\Delta_s+\varepsilon}.
\end{equation}

\noindent Following the various works of Vaughan, Br\"{u}dern and Wooley (see \cite{VaughanWooley1993,VaughanWooley1994,VaughanWooley1995,VaughanWooley2000,BruedernWooley2000,BrudernWooley2023Waringlarge}), we present some estimates regarding the mean values of smooth Weyl sums given by \eqref{DefUks}. To describe these estimates, we recall the concept of an \textit{admissible exponent} from the theory of smooth Weyl sums. A real number $\Delta_{k,s}$ is referred to as an \textit{admissible exponent} (for $k$) if it has the property that, whenever
$\varepsilon >0$ and $\eta = \eta(k,s,\varepsilon)$ is a sufficiently small positive number, then for any $2 \le R \le P^\eta$ where $P$ is sufficiently large, one has the bound

\begin{equation}
\label{Uks}
    U_{k,s}(P) \ll P^{s-k+\Delta_s+\varepsilon}.
\end{equation}

Also, if we use $\Delta_{k,s}^* = \Delta_{k,s}+10^{-10}$ instead of $\Delta_{k,s}$, then there exists an absolute constant $\eta >0$ such that whenever $R \le P^\eta$, one has

\begin{equation}\label{MVEeps}
\int_0^1 |g_k(\a,P,R)|^{s} \text{d}\a \ll P^{s-k+\Delta_{k,s}^*}.\end{equation}

\noindent We can therefore disregard $\varepsilon$ whenever dealing with such estimates, as the numerical data for $\Delta_{k,s}$ will be rounded up to the $6
^\text{th}$ decimal\footnote{\label{NoteDecimal}Rounded up to the $6^\th$ decimal as done in \cite{BruedernWooley2000}. Therefore, we align the computation to this order of magnitude, which explains the choice of $\Delta_{k,s}^*$ in \eqref{MVEeps}.}, and the result thus secured will depend on those data. Therefore, we shall not distinguish $\Delta_{k,s}$ from $\Delta_{k,s}^*$. Lastly, since our proofs contain only a finite number of steps, the real $\eta$ might be chosen at the very end so that whenever $R = P^\eta$ occurs, the results hold.

To ease future proofs, following \eqref{eq:defXk}, it turns out that upon choosing 

\begin{equation}
    \label{XkR}
    P = X_k = \floor{n^{1/k}} \et R= X_k^\eta,
\end{equation}

\noindent the previous estimate is of the shape

\begin{equation}
\label{MVEEven}
\int_0^1 |g_k(\a,X_k,R)|^{s} \text{d}\a \ll n^{\frac{s}{k}}n^{-1}n^{\frac{\Delta_{k,s}^*}{k}}.
\end{equation}

\noindent The term $n^\frac{s}{k}$ corresponds to the trivial estimate for the integrated sum. It is expected and proved (see, for instance, \cite{VaughanWooley1995}) that for $s$ sufficiently large in terms of $k$, there is no remaining term of shape $\frac{\Delta_{k,s}}{k}$. Hence, this quantity only recalls how far we are from the expected asymptotic.

The main literature produces data for $\Delta_{k,s}$ as long as $s$ is even (See \cite{VaughanWooley1995,VaughanWooley1993,VaughanWooley1994,VaughanWooley2000}). However, using Hölder's inequality, one can extend the computed values for $\Delta_{k,s}$ to any positive integer $s$ by convexity. Indeed, for any positive integer $s$ and real $\theta \in [0,1)$, let $S = 2(s+\theta) = \theta\cdot2(s+1)+(1-\theta)\cdot 2s$. Thus
\vspace{-5pt}
\begin{align}
\label{MVE}
    \int_0^1 |g_k(\a,X_k,R)|^{S} \text{d}\a &\ll \left(\int_0^1 |g_k(\a,X_k,R)|^{s} \text{d}\a\right)^{1-\theta}\left(\int_0^1 |g_k(\a,X_k,R)|^{s+1} \text{d}\a\right)^\theta \nonumber\\
    &\ll n^\frac{S}{k}n^{-1}n^\frac{\Delta_{k,S}}{k}
\end{align}
\leftcenters[0]{where}{$\Delta_{k,S} = \theta\Delta_{k,2(s+1)}+(1-\theta)\Delta_{k,2s}.$}

\noindent We present in Appendix \ref{AppAdmExp} admissible values for various $k$, $s$ even and the values needed for our computation. Also, note that some work done by Br\"{u}dern and Wooley \cite{BruedernWooley2000,Wooley2015} led them to break this convexity barrier. In particular, some of these novel values will be used in Sections \ref{FirstMajArc} and \ref{UAP}. Finally, from now on, we will write $g_k(\a)$ for $g_k(\a,X_k,R)$ with $X_k \et R$ fixed as in \eqref{XkR}.

As we just remarked, the mean value estimate is in a sense optimal when the exponent chosen is even. Therefore, if we go back to our integral problem over minor arcs, $f_2$ will be bounded using \eqref{EstF2Minor}, from which we will be left with

$$
\int_{\m(Q)}|g_4(\a,X_k,R)\cdots g_{158}(\a,X_k,R)|\text{d}\a.
$$

\noindent To bound such an integral, the secret lies in Cauchy's inequality and a suitable choice of $K \sqcup L =\set{4,\dots,158} $ to get
\begin{equation}
\label{EqExample}
    \int_{\m(Q)}|g_4(\a)\cdots g_{158}(\a)|\text{d}\a \ll \left(\int_{0}^1|G_K(\a)|^2\text{d}\a\right)^\half\left(\int_{0}^1|G_L(\a)|^2\text{d}\a\right)^\half,
\end{equation}

\leftcenters[0]{where}{$\displaystyle G_K = \prod\limits_{k\in K}g_k \et G_L = \prod\limits_{l \in L} g_l.$}

\leftcenters[0]{Let}{$\displaystyle \Kappa := \int_{0}^1|G_K(\a)|^2\text{d}\a$.}

\noindent It is therefore easy to perform a Hölder inequality of convexity with exponents $\{s_k\}_{k\in K}$ to reach that

$$
\Kappa \ll \prod\limits_{k \in K} U_{k,2s_k}(n^\frac{1}{k})^{\frac{1}{s_k}}.
$$

Then, using \eqref{MVEEven} whenever $s_k$ is even and \eqref{MVE} otherwise, Appendix \ref{AppAdmExp} provides admissible exponents so that
\begin{equation}
        \label{PrealphaK1}\Kappa \ll \prod\limits_{k \in K} n^{\frac{1}{k}\cdot \frac{2s_k-k+\Delta_{k,2s_k}}{s_k}} = \prod\limits_{k\in K}n^{\frac{2}{k}-\frac{1}{s_k}+\frac{\Delta_{k,2s_k}}{ks_k}} = n^{\kappa}n^{-1+\sum\limits_{k \in K}\frac{\Delta_{k,2s_k}}{ks_k}},
    \end{equation}

\leftcenters[0]{where}{$ \displaystyle
\kappa = \sum\limits_{k \in K}\frac{2}{k}.
$}

\noindent For reasons we shall see in the near future, we introduce the following quantity
\begin{equation}
\label{DefAlphaK}
\a(K) = 1 - \sum\limits_{k \in K}\frac{\Delta_{k,2s_k}}{ks_k}.
\end{equation}
We emphasise that \eqref{PrealphaK1} is better than the trivial estimate
$$
\Kappa \ll n^\kappa = G_K(0)^2=\prod\limits_{k\in K} g_k(0)^2,
$$
and leads to
\begin{equation}
    \label{alphaK1}
    \int_0^1|G_K(\a)|^2\text{d}\a \ll G_K(0)^2n^{-\a(K)}.
\end{equation}

Finally, we just showed that the bounding of a mixed integral as in \eqref{EqExample} relies on finding suitable sets $K \sqcup L = \set{4,\dots,158}$ to apply Cauchy's inequality and then to find suitable sets of exponents $\set{s_k}_k$ to use Hölder's inequality. Ford's algorithm (see \cite{Ford1995UnlikePowersI,Ford1996UnlikePowersII}) and a revisited version of the code of \cite{KLX2022} are precisely designed to determine such a dissection and exponents.
As $4$ and $10$ play a pivotal role in our proof and will not be placed in the same subset, we denote by $\set{4,\dots,158} = K_4 \sqcup K_{10}$ any disjoint union where $K_4$ contains $4$ and $K_{10}$ contains ${10}$.
From now on, for any such $K_4 \et K_{10}$, we write
$$
\cf(\a) = f_2(\a)G_{4}(\a)G_{10}(\a),
$$

\leftcenters[0]{where}{$
G_4 = \prod\limits_{k\in K_4}g_k \et G_{10} = \prod\limits_{k\in K_{10}}g_k.
$}

\noindent Using Ford's method, we take

\begin{align}
    \label{DefKF}
    K_{4} := &\{4,6,8,22,30,34,40,50,52,60,62,66,68,70,74\}\cup\{88,\dots,96\} \nonumber \\
    &\cup \{100,106,108,110,118,120,130,132,134,144,152,154,156\},
\end{align}

\noindent and then introduce its complementary $K_{10}$ as

\begin{align}
    \label{DefKS}
    K_{10} = &\set{10,\dots,20} \cup \{24,\dots,32 \} \cup \{36,38 \} \cup\{42,\dots,48 \}\cup \{54,56,58,64,72 \}\cup \{76,\dots, 86 \} \nonumber \\
    &\cup \{ 98,102,104,112,114,116\} \cup \{122,\dots, 128 \} \cup \{136,\dots,142 \} \cup \{146,148,150,158 \}.
\end{align}

\noindent Therefore, we obtain the first minor arc estimate of our process, whose proof will be presented in the next section. 

\begin{proposition} \label{LemmaMajArcFirst}
Let $\delta = 4.3\times 10^{-5}$ and  $Q_1 = n^{1/2-\delta/2}$. Then we have
    $$\int_{\m(Q_1)} |f_2(\a)G_4(\a)G_{10}(\a)| \text{d}\a \ll \cf(0)n^{-1-\delta/4}.$$
\end{proposition}

\noindent To prove \autoref{LemmaMajArcFirst}, we need the following lemma

\begin{lem}
\label{LemmaK4&K6}
    Let $\delta_4 = 3.3\times10^{-5}$. Then one has

    $$\int_0^1 |G_{4}(\a)|^2 \text{d}\a \ll n^{-0.76-\delta_4}G_4(0)^2.$$

    \noindent Similarly, let $\delta_{10} = 1.0\times10^{-5}$. Then, one has
    
    $$\int_0^1 |G_6(\a)|^2 \text{d}\a \ll n^{-0.74-\delta_{10}}G_{10}(0)^2.$$
\end{lem}

\section{First candidate for the minor arc bound}
\label{FirstMajArc}

Following the previous section, we define $$\Kappa_j = \int_0^1 |G_j(\a)|^2 \text{d}\a \quad (j = 4 \ou 10).$$
Similarly, for $j = 4 \ou 10$, let $\kappa_j = \sum\limits_{k\in K_j} \frac{2}{k}.$ Then \eqref{alphaK1} yields

$$
\Kappa_4 \ll G_4(0)^2n^{-\a(K_4)},
$$
\leftcenters[0]{and}{$\Kappa_{10} \ll G_{10}(0)^2n^{-\a(K_{10})}.$}

\noindent Thus, it remains to prove that $\a(K_4) \et \a(K_{10})$ are bigger than the exponents announced in \autoref{LemmaK4&K6}.

\begin{proof}{\textit{(Of \autoref{LemmaK4&K6})}}
    As we mentioned in \eqref{DefAlphaK} and \eqref{alphaK1}, we use the set $\set{s_k}_{k\in K_4}$ as presented in \autoref{SetK4}

\begin{table}[H]
\centering

\begin{minipage}{0.23\textwidth}
\centering
\begin{tabular}{|c|c|}
\hline
$k$ & $s_k$\\
\hline
4 & $s_4$\\
6 & 6\\
8 & 7\\
22 & 22\\
30 & 30\\
34 & 35\\
40 & 40\\
50 & 50\\
52 & 52\\
\hline
\end{tabular}
\end{minipage}
\hspace{-25pt}
\begin{minipage}{0.23\textwidth}
\centering
\begin{tabular}{|c|c|}
\hline
$k$ & $s_k$\\
\hline
    60 & 60\\
    62 & 62\\
    66 & 67\\
    68 & 70\\
    70 & 72\\
    74 & 76\\
    88 & 90\\
    90 & 92\\
\hline
\end{tabular}
\end{minipage}
\hspace{-25pt}
\begin{minipage}{0.23\textwidth}
\centering
\begin{tabular}{|c|c|}
\hline
$k$ & $s_k$\\
\hline
    92 & 93 \\
    94 & 95 \\
    96 & 97\\
    100 & 102\\
    106 & 107\\
    108 & 110\\
    110 & 111\\
    118 & 120\\
\hline
\end{tabular}
\end{minipage}
\hspace{-25pt}
\begin{minipage}{0.23\textwidth}
\centering
\begin{tabular}{|c|c|}
\hline
$k$ & $s_k$\\
\hline
    120 & 122\\
    130 & 132 \\
    132 & 134\\
    134 & 136\\
    144 & 146\\
    152 & 153\\
    154 & 156\\
    156 & 158\\
\hline
\end{tabular}
\end{minipage}

\caption{Hölder's exponents for $K_4$}
\label{SetK4}
\end{table}
\noindent where $s_{4} >0$ is chosen such that
\begin{equation*}
\sum\limits_{k \in K_4}\frac{1}{s_k}=1.     
\end{equation*}
We note that $s_{4} = 3.650062\dots$ Using admissible exponents as presented in Table \ref{TabAdmExpK4} in the Appendix, we~have
\begin{equation}
    \label{Kappa4}
    \Kappa_4 \ll n^{\kappa_4-\alpha(K_4)},
\end{equation}
where $\alpha(K_4) = 0.760033\dots > 0.76 + \delta_4.$

Regarding $\Kappa_{10}$, we use the set of exponents introduced in \autoref{SetK10} with $s_{10} = 9.991737 \dots$ given such that
\vspace{-15pt}

$$\sum\limits_{k \in K_{10}}\frac{1}{s_k}=1,$$

\noindent along the admissible exponents of Table \ref{TabAdmExpK10} to reach
\begin{equation}
    \label{Kappa10}
    \Kappa_{10} \ll n^{\kappa_{10}-\alpha(K_{10})},
\end{equation}
where $\alpha(K_{10}) = 0.7400103\dots > 0.74 + \delta_{10}.$
This completes the proof.

\begin{table}[H]
\centering

\begin{minipage}{0.23\textwidth}
\centering
\begin{tabular}{|c|c|}
\hline
$k$ & $s_k$\\
\hline
10 & $\mathbf{s_{10}}$\\
12 & 12\\
14 & 14\\
16 & 16\\
18 & 18\\
20 & 20\\
24 & 24\\
26 & 26\\
28 & 29\\
\hline
\end{tabular}
\end{minipage}
\hspace{-25pt}
\begin{minipage}{0.23\textwidth}
\centering
\begin{tabular}{|c|c|}
\hline
$k$ & $s_k$\\
\hline
32 & 32.5\\
36 & 37.5\\
38 & 39\\
42 & 43\\
44 & 45.5\\
46 & 47.25\\
48 & 49\\
54 & 55\\
56 & 58\\
\hline
\end{tabular}
\end{minipage}
\hspace{-25pt}
\begin{minipage}{0.23\textwidth}
\centering
\begin{tabular}{|c|c|}
\hline
$k$ & $s_k$\\
\hline
58 & 60\\
64 & 66\\
72 & 74\\
76 & 78\\
78 & 80\\
80 & 83\\
82 & 84\\
84 & 86\\
86 & 88\\
\hline
\end{tabular}
\end{minipage}
\hspace{-25pt}
\begin{minipage}{0.23\textwidth}
\centering
\begin{tabular}{|c|c|}
\hline
$k$ & $s_k$\\
\hline
98 & 101\\
102 & 104\\
104 & 106\\
112 & 115\\
114 & 117\\
116 & 120\\
122 & 125\\
124 & 127\\
126 & 129\\
\hline
\end{tabular}
\end{minipage}
\hspace{-25pt}
\begin{minipage}{0.23\textwidth}
\centering
\begin{tabular}{|c|c|}
\hline
$k$ & $s_k$\\
\hline
128 & 131\\
136 & 139\\
138 & 143\\
140 & 145\\
142 & 147\\
146 & 150\\
148 & 151\\
150 & 155\\
158 & 162\\
\hline
\end{tabular}
\end{minipage}
\caption{Hölder's exponents for $K_{10}$}
\label{SetK10}
\end{table}
\end{proof}

\begin{proof}{(\textit{\autoref{LemmaMajArcFirst}})}
Let $1 \le Q \le n^\half$. Using Cauchy's inequality, one has

$$
\int_{\m(Q)}  |\cf(\a)|\text{d}\a \ll \sup\limits_{\a \in \m(Q)} |f_2(\a)|\left(\int_0^1 |G_4(\a)|^2 \text{d}\a\right)^\half\left(\int_0^1|G_{10}(\a)|^2\text{d}\a\right)^\half.
$$

\noindent The first expression on the right is estimated by \eqref{EstF2Minor} whereas Lemma \ref{LemmaK4&K6} provides bounds for the two remaining integrals. Therefore, we get

$$
\int_{\m(Q)} |\cf(\a)|\text{d}\a \ll \cf(0) Q^{-\half} n^{-\frac{3}{4}-\frac{\delta_{10}+\delta_4}{2}+\varepsilon}.
$$

\noindent Choosing $Q = Q_1 = n^{1/2-\delta/2}$, and because $\delta_{10}+\delta_4 = \delta$, we get

$$
\int_{\m(Q_1)} |\cf(\a)|\text{d}\a \ll \cf(0)n^{-1-\frac{\delta}{4}+\varepsilon}.
$$

\noindent This ends the proof of Lemma \ref{LemmaMajArcFirst}.
    
\end{proof}

We emphasise that $\delta$ has the following property

\begin{equation}
    \label{delta}
    \delta < \a(K_4)+\a(K_{10}) - \frac{3}{2}.
\end{equation}

\section{Mean value estimate over intermediate arcs}
\label{Pruning}

We start with some recent technological advances in the domain of smooth Weyl sum estimates. First, we restate the following crucial result achieved by Br\"{u}dern and Wooley \cite[Th 4.2]{BrudernWooley2023Waringlarge}.

\begin{theorem}
\label{Theorem 4.2}
Suppose that $s$ is a real number with $s\ge 2$ and $\Delta_{k,s}$ is an admissible
exponent. \\ Then whenever $Q$ is a real number with $1 \le Q \le P^{k/2}$, one has the uniform bound

$$
V_{k,s}(Q;P,R) :=\int_{\M(Q,P)}|g_k(\a;P,R)|^s \text{d}\a \ll P^{s-k}Q^{2\Delta_{k,s}/k+\varepsilon}.
$$

\end{theorem}

The proof of Lemma \ref{LemmaMajArcLog} is now within our grasp, up to replacing $f_2$ by its completed version $f_2^*$. To this aim, we introduce
\vspace{-3pt}
$$E_2(\a) := f_2(\a)-f_2^*(\a),$$

\noindent for any $\a \in [0,1)$. We make use of \cite[Theorem 4.1]{vaughan1981hardy}, which gives that
\vspace{-1pt}
\begin{align}
\label{Ek}
    E_2(\a) &\ll q^{\half + \varepsilon}(1+n\b)^\half, \nonumber \\
    \text{and even } E_2(\a) &\ll q^{\half + \varepsilon} \text{ when } |\b| \ll \frac{n^{\half\minus}}{2kq}.
\end{align}

\noindent Hence, whenever $\a$ belongs to $\M(Q_1)$, \eqref{Ek} yields
\vspace{-1pt}
\begin{equation}\label{EstimateDelta2}
E_2(\a) \ll q^{\half+\varepsilon} \ll Q_1^{\half + \varepsilon}.
\end{equation}

\noindent Using this estimate, one can easily replace $f_2$ by $f_2^*$.

\begin{lem}
\label{LemmaReplacef2}
    One has
    $$
    \int_{\M(Q_1)} |E_2(\a)G_4(\a)G_{10}(\a)|\text{d}\a \ll \cf(0) n^{-1-\frac{\delta_1}{2}}.
    $$
\end{lem}

\begin{proof}
    The proof of this Lemma is very similar to that of Lemma \ref{LemmaMajArcFirst}. Changing $f_2$ by $E_2$ in the proof and using \eqref{EstimateDelta2} in place of \eqref{EstF2Minor}, it is straightforward to obtain
    \begin{align*}
        \int_{\M(Q_1)} |E_2(\a)G_4(\a)G_6(\a)|\text{d}\a &\ll Q_1^{\half + \varepsilon}G_4(0)G_6(0) n^{-\frac{3}{4}-\delta/2} \\
        &\ll \cf(0) n^{-\half}n^{\half(\half-\frac{\delta}{2})+\varepsilon}n^{-\frac{3}{4}-\delta/2}\\
        & \ll \cf(0) n^{-1-\frac{3\delta}{4}+\varepsilon}.
    \end{align*}
\end{proof}

Now, we introduce a mixed mean value that will be of great interest as we expect to study integration over $\M(Q_1)\backslash\M(X)$ for some smaller $X$ to prune the major arcs. Let $1 \le Q \le Q_1$, then define, for $j = 4 \ou {10}$, the quantity

\begin{equation}
\label{UpsX}
    \Upsilon_j(Q_1,Q) := \int_{\M(Q_1)\backslash\M(Q)} |f_2^*(\a)G_j^2(\a)|\text{d}\a.
\end{equation}

\noindent The core idea of the pruning step is based on a smart dyadic argument applied to the intermediate arc on which $\Upsilon_j(Q_1,Q)$ is defined. We write $\Upsilon_j(Q) := \Upsilon_j(Q,Q/2)$. 

\begin{lemma}
\label{UpsNQ}
Let $1 \le Q \le Q_1$. One has
$$
\int_{\cn(Q)} |f_2^*(\a)G_4(\a)G_{10}(\a)|\text{d}\a \ll \cf(0)n^\minus Q^{-\delta}.
$$
\end{lemma}

\begin{proof}
    Using Cauchy's inequality, we have
    
    \begin{equation}
    \label{CSUps}
    \int_{\cn(Q)} |f_2^*(\a)G_4(\a)G_{10}(\a)|\text{d}\a \ll \Upsilon_4(Q)^\half\Upsilon_{10}(Q)^\half.
    \end{equation}
    
    \noindent Thus, it is sufficient to prove that for any $ 1 \le Q \le Q_1$ and both $j = 4 $ and $10$, we have
    
    \begin{equation}
    \label{UpsjNQ}
    \Upsilon_j(Q) \ll f_2(0)G_j(0)^2n^{-1}Q^{2(1-\a(K_j))-\half+\varepsilon},
    \end{equation}
    
    \noindent for any $1 \le Q \le Q_1$ and for $j = 4 \ou 6$.

    Now, fix $\a$ in $\cn(Q) \subset \m(Q)$. In view of the definition of $f_2^*$ given in \eqref{DefFk*}, we use \eqref{EstSk} and \eqref{Estnuk} to reach
    
    $$
    f_2^*(\a) \ll \frac{1}{q}S_2(q,a)\nu_2(\b) \ll q^{-\half+\varepsilon}\frac{X_2}{(1+|n\b|)^\half}.
    $$
    
\noindent Now, as $f_2(0) \asymp X_2$ and $\a \in \cn(Q)$, using \eqref{EstNQ}, one can conclude that
    
    \begin{equation}
        \label{EstF2*}
        \sup\limits_{\a \in \cn(Q)}|f_2^*(\a)| \ll f_2(0)\frac{Q^{\varepsilon}}{q^\half(1+n|\b|)^\half}\ll f_2(0)Q^{-\half+\varepsilon}.
    \end{equation}

We prove the result for $j = 4$, the exact same argument will apply for $j = 10$. Using \eqref{EstF2*} and Hölder's inequality with coefficients $\set{s_k, k \in K_j}$ defined in \autoref{SetK4}, one has

    \begin{align*}
        \Upsilon_4(Q) &\ll f_2(0)Q^{-\half+\varepsilon}\int_{\M(Q)}\prod\limits_{k \in K_4}|G_4(\a)|^2\text{d}\a \\
        &\ll f_2(0)Q^{-\half+\varepsilon}\prod\limits_{k \in K_4}V_{k,2s_k}(Q)^\frac{1}{s_k} \\
        &\ll f_2(0)Q^{-\half+\varepsilon} \prod\limits_{k\in K_4} \left(X_k^{2s_k}n^{-1}Q^{\frac{2\Delta_{2s_k}}{k}+\varepsilon}\right)^\frac{1}{s_k} \quad \quad \text{(\autoref{Theorem 4.2})} \\
        &\ll f_2(0)G_4(0)^2n^{-1}Q^{\sum\limits_{k\in K_4}\frac{2\Delta_{2s_k}}{ks_k}-\half+\varepsilon}.
    \end{align*}
    
\noindent Finally, \eqref{DefAlphaK} determines the value of $\sum\limits_{k\in K_4}\frac{\Delta_{2s_k}}{ks_k}$ in terms of $\a(K_4)$. Then one easily concludes that

    $$\Upsilon_4(Q) \ll f_2(0)G_4(0)^2n^{-1}Q^{2(1-\a(K_4))-\half}.$$
    
\noindent Implementing both bounds of \eqref{UpsjNQ} and using the definition of $\delta$ given by \eqref{delta} into \eqref{CSUps}, we obtain
\vspace{-1pt}
    \begin{align*}
        \int_{\cn(Q)} |f_2^*(\a)G_4(\a)G_{10}(\a)|\text{d}\a &\ll \cf(0)n^\minus Q^{(2-\a(K_4)-\a(K_{10}))-\half} \\
        &\ll \cf(0)n^\minus Q^{-\delta}.
    \end{align*}
\end{proof}

We are now able to shrink the integration domain to $\M(Z)$ where $Z = L^{1/15}$.

\begin{lemma}
    \label{Lemmashrink}
    One has
    $$
    \int_{\M(Q_1)\backslash\M(Z)}|f_2^*(\a)G_4(\a)G_{10}(\a)|\text{d}\a \ll \cf(0)n^{-1}L^{-\rho}
    $$
    for some positive $\rho$.
\end{lemma}

\begin{proof}
Note that $$\M(Q_1)\cap\m(Z) \subset\bigcup\limits_{j = 0}^J \cn(2^{-j}Q_1),$$ where $J = \floor{\frac{\log Q_1/Z}{\log 2}} \ll \log n$. Therefore, using \autoref{UpsNQ}, this dyadic dissection leads to
\begin{align*}
\int_{\M(Q_1) \cap \m(Z)} |f_2^*(\a)G_4(\a)G_{10}(\a)|\text{d}\a 
&\le \sum\limits_{j=0}^J \int_{\cn(2^{-j}Q_1)} |f_2^*(\a)G_4(\a)G_{10}(\a)|\text{d}\a\\
&\ll\cf(0)n^\minus \sum\limits_{j=0}^J(2^{-j}Q_1)^{-\delta} \\
&\ll \cf(0)n^\minus Z^{-\delta}.    
\end{align*}

\noindent Therefore, as $Z =L^\frac{1}{15}$, one gets the required result.
    \end{proof}
Combining Lemmata \ref{LemmaMajArcFirst}, \ref{LemmaReplacef2} and \ref{Lemmashrink}, we thus prove \autoref{LemmaMajArcLog}, that is
$$T(n;\eta) = \int_{\M(L^{1/15})} f_2^*(\a)G_4(\a)G_{10}(\a)\text{d}\a +\bigoh{\cf(0)n^\minus L^{-\rho}},$$
\noindent where $\rho = \delta/15$.

\section{Treatment of the major arcs}
\label{MajArcsasymptotic}
For a more convenient approach, we take the set of major arcs $\Kfrak$ to be the union of the intervals
$$
\Kfrak(q,a) = \set{\a \in [0,1) : |\a - \frac{a}{q}| \le \frac{L^\frac{1}{15}}{n}},
$$
with $0 \le a \le q \le L^\frac{1}{15} \et (a,q)=1$. The set of minor arcs complementary to $\Kfrak$ is then $\kfrak = [0,1) \backslash \Kfrak$. \\ We write
$$
\ncurly(n,\eta) = \int_{\Kfrak} f_2^*(\a)G_4(\a)G_{10}(\a) e(-n\a)\text{d}\a.
$$
So far, using Lemma \ref{LemmaMajArcLog} and remarking that $\m(L^\frac{1}{15}) \subset \kfrak$, we have
\begin{align}
    \label{finalmajorarc}
\int_0^1 \cf(\a) e(-n\a)\text{d}\a &= \ncurly(n,\eta) + \bigoh{\cf(0)n^{-1}L^{-\rho}}.
\end{align}
In the same spirit as what we already did to replace $f_2$ by its complete sum, the estimate we reach for the associated $E_k$ is non-trivial only for small values of $Q$, which is why we reach such a convenient major arc set. The standard theory of smooth Weyl sums (see \cite[Lemma 5.4]{Vaughan1989} or \cite[Lemma 8.5]{Wooley1991} for a similar conclusion) will take care of the remaining smooth exponential sums as follows. 

\begin{lemma}
\label{SmoothWeilCousin}
There is a positive number $c = c(\eta)$ such that for $4 \le k \le 158$, and $\a \in \Kfrak (q,a)\subset \Kfrak$, one has
$$
g_k(\a) = cq^\minus S_k(q,a)\nu_k(\b) + \bigoh{X_kL^{-1/4}},
$$
where the error term is uniform in $\a$ and $S_k \et \nu_k$ are defined in \eqref{DefFk*}.
\end{lemma}
Therefore, for $k \ge 4$, we introduce
\begin{equation}
    g_k^*(\a,q,a) = c q^\minus S_k(q,a)\nu_k(\b).
\end{equation}
For brevity, we write
\begin{equation}
\label{gkstar}
    g_k^*(\a) = \left\{ \begin{split}
        g_k^*(\a,q,a) \quad &\text{ for } \a \in \Kfrak,\\
        0 \quad \quad \quad \quad \quad &\text{ else},
    \end{split}\right.
\end{equation}
and
\begin{align}
\label{Ekg}
    E_k(\a) = g_k(\a)-g_k^*(\a) \ll X_kL^{-1/4}.
\end{align}
Also note that \eqref{EstSk} and \eqref{Estnuk} yield the following result.
\begin{lemma}{\textit{\cite[Lemma 5.4]{vaughan1981hardy}}}
    \label{EstWk}
    For $k \ge 2, 1 \le Q \le n^\half$ and $\a \in \M(Q)$, we have 
    $$
    g_k^*(\a) \ll \left(\frac{n}{q}\right)^\frac{1}{k}(1+n|\b|)^{-\frac{1}{k}}.
    $$
\end{lemma}
We take full advantage of \autoref{SmoothWeilCousin} to replace $g_k$ by its related function $g_k^*$. 
To perform such a replacement, we proceed to change them one by one.
\begin{lemma}
\label{ReplaceSmoothCousin}
    One has
    $$
    \ncurly(n,\eta) = \int_{\Kfrak}f_2^*(\a)G_4^*(\a)G_{10}^*(\a)e(-n\a) \text{d}\a + \bigoh{\cf(0)n^\minus L^{-\frac{1}{20}}},
    $$
    where $G_j^* = \prod\limits_{k \in K_j}g_k^*$ for $j = 4 \ou {10}$.
\end{lemma}

\begin{proof}

Let $\a \in \Kfrak$ and $k \ge 4$. \autoref{EstWk} leads to 

$$g_k^*(\a) \ll n^\frac{1}{k}.$$

\noindent Therefore, as $\Kfrak$ is of Lebesgue measure $n^\minus L^{\frac{1}{5}}$, substituting the trivial estimate for $f_2^*$, the trivial estimates for $g_k \et g_k^*$ when needed and using \eqref{Ekg}, we obtain
\vspace{-1pt}
\begin{align*}
    \int_{\Kfrak} |f_2^*(\a)g_4^*(\a)\cdots g_{2(k-1)}^*(\a)E_{2k}(\a)g_{2(k+1)}(\a)\cdots g_{158}(\a)| \text{d}\a 
    &\ll \cf(0)n^{-\frac{1}{2k}} \int_{\Kfrak} E_{2k}(\a)\text{d}\a \\
    & \ll \cf(0) L^{-\frac{1}{4}}n^\minus L^\frac{1}{5} \\
    &\ll \cf(0) n^\minus L^{-\frac{1}{20}}.
\end{align*}

\noindent Therefore, replacing one by one the $g_k$'s, we save a positive power of $\log n$ at each step. The result is secured in 78 such replacements, since
\vspace{-1pt}
\begin{align*}
    \Big|\ncurly(n,\eta) - \int_{\Kfrak}&f_2^*(\a)G_4^*(\a)G_6^*(\a)e(-n\a) \text{d}\a \Big| \\
    &\ll \sum\limits_{k=2}^{79} \int_{\Kfrak} |f_2^*(\a)g_4^*(\a)\cdots g_{2(k-1)}^*(\a)E_{2k}(\a)g_{2(k+1)}(\a)\cdots g_{158}(\a)| \text{d}\a \\[2mm]
    &\ll \cf(0) n^\minus L^{-\frac{1}{20}}.
\end{align*}

\end{proof}

At this stage of the proof, the final step is just routine. Introduce 

$$
A(n,q) = \frac{c^{78}}{q^{79}}\sum\limits_{(a,q)=1} S_2\cdots S_{158}(q,a)e(-an/q),
$$

\noindent where $c$ has been defined in \autoref{SmoothWeilCousin}
and

$$
\S(n,Y) = \sum\limits_{q \le Y}A(n,q). 
$$

\leftcenters[0]{Let also}{$\I(n,Y) = \int_{-Y/n}^{Y/n} \nu_2(\b)\cdots \nu_{158}(\b)e(-n\b)\text{d}\b .$}

\noindent The next step is to replace the integral over the major arc by an integral over $[0,1)$ as the various $g_k^*$ will be small outside of $\Kfrak(Y,q,a)$ for $1 \le Y \le \frac{n^\half}{4}.$The following arguments are similar to those of \cite[Equations (4.33) and (4.34)]{Ford1996UnlikePowersII}, but we reproduce them for the sake of completeness. Define the quantity

\begin{equation}
    \label{mu}
    \mu = \sum\limits_{k =1}^{79}\frac{1}{2k} = 2.48\dots
\end{equation}

\noindent The most pedestrian way to introduce $\S(n,Z) \et \I(n,Z)$ is described below. Regarding the study of the integral over $\Kfrak$, we see that it can be decomposed into a sum of smaller integrals over $\Kfrak(q,a)$. Using \eqref{gkstar}, we obtain
\begin{align*}
    \int_{\Kfrak} f_2^*(\a)g_4^*(\a)\cdots g_{158}^*(\a) \text{d}\a 
    &=\sum\limits_{q \le Z}\sum\limits_{(a,q) =1} \int_{\Kfrak(q,a)} f_2^*g_4^*\cdots g_{158}^*(\a,q,a) \text{d}\a  \\
    &= \sum\limits_{q \le Z}\sum\limits_{(a,q) =1} \int_{-Z/n}^{Z/n} \frac{c^{78}}{q^{79}}\prod\limits_{k=1}^{79}S_{2k}(q,a)\prod\limits_{k=1}^{79} \nu_{2k}(\b) e\left(-n\b-n\frac{a}{q}\right)\text{d}\b \\
    &= \S(n,Z)\I(n,Z).
\end{align*}

\noindent Then, from the definition of $g_k^*$ and $f_2^*$, we define $\ncurly^*$ so that
\begin{equation}
\label{SnIn}
    \ncurly^*(n,Z) := \sum\limits_{q \le Z}\sum\limits_{(a,q) =1}\int_0^1 f_2^*(\a,q,a)G_4^*(\a)G_{10}^*(\a)e(-n\a)\text{d}\a = \S(n,Z)\I\left(n,Z\right).
\end{equation}

\noindent The final steps are to replace both the singular integral and the singular series by their full version and prove that they are strictly positive.
\subsection{Singular series}
\label{SecSingularSeries}

It is well known that for any positive integer $n$, the map $q \mapsto A(n,q)$ is a multiplicative function.

\begin{lemma}
\label{Singular series}

The truncated singular series $\S(n,Y)$ converges to $\S(n) = \S(n,\infty)$ as $Y $ goes to $+\infty$, at the rate
\vspace{-1pt}
$$
\S(n,Y) = \S(n) + \bigoh{Y^{2-\mu}},
$$

\noindent where $\mu$ is as in \eqref{mu}. Also,
\vspace{-1pt}
$$
\S(n) \gg 1.
$$
\end{lemma}

\begin{proof}

First, let $p$ be a prime, and let $l \ge 1$ be an integer. Then, using \eqref{EstSk}, we see that
\vspace{-1pt}
\begin{align*}
A(n,p^l) 
&= \frac{1}{(p^l)^{79}}\sum\limits_{\substack{1 \le a \le p^l \\ (a,p^l)=1}}\prod\limits_{k = 1}^{79} S_{2k}(p^l,a)e(-an/p^l) \\
&\ll \frac{1}{(p^l)^{79}}p^l\prod\limits_{k=1}^{79} (p^l)^{1-\frac{1}{2k}+\varepsilon} \\
&\ll (p^l)^{1-\mu+\varepsilon}.
\end{align*}

\noindent Hence, by multiplicativity,

\begin{equation}
    \label{MajAnq}    
A(n,q) \ll q^{1-\mu+\varepsilon},
\end{equation}

\noindent so that for $\varepsilon$ not too large, $\S(n,Y)$ converges, say, to $\S(n)$. In particular, \eqref{MajAnq} shows,

$$
\S(n,Y) = \S(n) + \bigoh{Y^{2-\mu+\varepsilon}}.
$$

\noindent Additionally, by multiplicativity, one can write $\S(n)$ as a Euler product

$$
\S(n) = \prod\limits_{k = 1}^{+\infty} \chi_p,
$$

\noindent where $\chi_p = \sum\limits_{l=0}^{+\infty} A(n,p^l)$.
Local solubility is resolved using \cite[Lemma 6.4]{Ford1995UnlikePowersI}. For each prime $p$ and integer $k$, let $p^\tau || \ k$ and define
\vspace{-1pt}
$$
\gamma(k,p) = \left\{\begin{split}
    \tau +2 &\si p=2 \et \tau >0, \\
    \tau +1 &\text{ otherwise}.
\end{split}\right.
$$

\noindent Let $k_1 = 2, \dots, k_{79} = 158$. Therefore, we distinguish between the following three cases:

\leftcenters[0]{We start when $p \ge 5$, then}{$
\sum\limits_{i = 1}^r \frac{p-1}{p^{\gamma(k_i,p)}(k_i,p-1)} \ge  \frac{p-1}{p} \sum\limits_{i\in\{1,2,3,4,6,8,9,12\}}\frac{1}{2i} \ge 1.
$}

\leftcenters[0]{Also, for {$p = 3$}, we have}{$
\sum\limits_{i = 1}^r \frac{3-1}{3^{\gamma(k_i,3)}(k_i,3-1)} \ge \frac{2}{3}\sum\limits_{i \in \{1,2,4\}} \frac{1}{2} \ge 1.
$}

\leftcenters[0]{Lastly, {$p = 2$} yields}{$\sum\limits_{i = 1}^r \frac{2-1}{2^{\gamma(k_i,2)}(k_i,2-1)} \ge \sum\limits_{i \in \{1,3,5,7,9,11,13,17\}} \frac{1}{2^3} \ge 1.
$}

\noindent Thus, with $r = 79 \et \gamma = 9$ (attained for $p=2 \et k = 2^7 = 128$), then for any $p$ prime, \cite[Lemma 6.4]{Ford1995UnlikePowersI} gives that $\chi_p$ is real, giving that $\S(n)$ is real and $\chi_p$ also satisfies

\begin{equation}
    \label{chipPetit}
    \chi_p \gg p^{-711}.
\end{equation}

\noindent Finally, note that for any $p$ prime, \eqref{MajAnq} ensures that

\begin{equation}
    \label{chipGrand}
    |1-\chi_p| \le \sum\limits_{h=1}^{+\infty} |A(n,p^h)| \le p^{1-\mu},
\end{equation}

\noindent which means that the infinite product $\prod\limits_{p\in \P}\chi_p$ converges. Combining equations \eqref{chipPetit} and \eqref{chipGrand}, there exist positive constants $c_0$ and $c_1$ such that

$$
\S(n) =  \prod\limits_{p > c_0} \chi_p \prod\limits_{p \le c_0} \chi_p\gg  \prod\limits_{p > c_0} (1-c_1p^{1-\mu}) \prod\limits_{p\le c_0}p^{-711} \gg 1.
$$
\end{proof}

\subsection{Singular integral}
\label{SecSingularIntegral}
Recall that $$\displaystyle \I(n,Z) = \int_{-Z/n}^{Z/n} \nu_2(\b)\cdots\nu_{158}(\b)e(-n\b)\text{d}\b.$$
In view of \eqref{Estnuk}, we get
\begin{align*}
    \I(n)-\I(n,Z) 
    &\ll \cf(0)\int_{Z/n < |\b| < \half }\frac{1}{(1+n|\b|)^\mu}\text{d}\b \\
    &\ll \cf(0)n^\minus Z^{1-\mu}.
\end{align*}
A similar approach following \cite[Theorem 2.3]{vaughan1981hardy} shows that
\begin{equation}
\label{SingularIntegral}
    \I(n) = \frac{\prod\limits_{k=1}^{79}\Gamma\left(1+\frac{1}{2k}\right)}{\Gamma(\mu)}\cf(0)n^\minus\left(1+\bigoh{n^{-\frac{1}{158}}}\right).
\end{equation}
Thus, the conclusions of Lemma \ref{Singular series}; equations \eqref{MajAnq} and \eqref{SingularIntegral} applied to \eqref{SnIn} lead to the existence of a positive number $\nu$ for which
\begin{align}
\label{10.7}
\ncurly^*(n,Z) 
&= \S(n,Z)\I(n,Z) \nonumber\\
&= \S(n)\I(n) + \bigoh{\cf(0)n^\minus L^{-\nu}}.
\end{align}
We finally conclude by \eqref{finalmajorarc} and \eqref{10.7} that
$$
\ncurly(n) = \S(n)\I(n) + \bigoh{\cf(0)n^\minus L^{-\tau}},
$$
for some strictly positive and sufficiently small $\tau$. Another remark is that \autoref{Singular series} and \eqref{SingularIntegral} ensure that the main term predominates. Since $T(n,\eta)$ diverges as $n$ tends to infinity, this implies that $T(n;\eta) \ge1$ for $n$ large enough. The proof of \autoref{ThEvenUnlikePowers} is now complete.

\section{Unlike ascending powers}
\label{UAP}

We aim to prove the qualitative improvement of \cite{LiuZhao2021_RepresentationUnlikePowers}.

\begin{theorem}
Let $n$ be a sufficiently large integer. Then for some $\eta$ sufficiently small, $n$ can be represented in the form
    \begin{equation}
    \label{EqUnlikePowers}
        x_1^2 + x_2^3 + \cdots + x_{11}^{12} + x_{12}^{13} + x_{13}^{17} = n,
    \end{equation}
    
  \noindent with $1 \le x_1 \le n^\half, 1 \le x_2 \le n^{\frac{1}{3}}$, $x_{i-1} \in \A(n^{1/i},n^{\eta/i})$ for $ 4 \le i \le 13$ and $x_{13} \in \A(n^{1/17},n^{\eta/17})$. In particular, if we denote by $H(n;\eta)$ the number of such solutions weighted as in \eqref{weight}, then
  \begin{equation}
  H(n;\eta) = c^{11}\S(n)\frac{\Gamma(1+1/2)\cdots \Gamma(1+1/13)\Gamma(1+1/17)}{\Gamma(\theta)}n^{\theta-1} + \bigoh{n^{\theta-1}(\log n)^{-\rho}}
  \end{equation}
  for some positive $\rho$, $c = c(\eta)$  and $\theta = \sum\limits_{k=2}^{13} \frac{1}{k} + \frac{1}{17} = 2.239\dots$
\end{theorem}

Regarding the structure of the proof, we arrive at this result by combining ideas drawn from \cite{LiuZhao2021_RepresentationUnlikePowers,BrudernWooley2023Freiman}, and we make use of some recent improvements concerning smooth Weyl sums, as in the previous sections. To this end, let $\a \in [0,1)$ be any real number, and define
\begin{equation}
    \label{weight}
F_3(\a) = \sum\limits_{1 \le x \le X_3} w(2x/X_3)e(\a x^3),
\end{equation}

\leftcenters[0]{where}{$\displaystyle w(t) = \exp\left(-\frac{1}{1-(t-1)^2}\right).$}

\leftcenters[0]{Let us write}{$\displaystyle \cf(\a) = F_2(\a)F_3(\a)g_4(\a)\cdots g_{13}(\a)g_{17}(\a).$}

\noindent Then, using tabulated values of Appendix C, we reach the first non-trivial major arc estimate.

\begin{lemma}
\label{EstMdelta}
Let $\delta = 12.782\times10^{-3}$ and define $Q_1 = n^{\half-\delta}$. Then one has
$$
\int_{\m(Q_1)} |\cf(\a)|\text{d}\a \ll \cf(0)n^\minus n^{-\delta+\varepsilon}.
$$
\end{lemma}

\noindent Lemma \ref{EstMdelta} will be proved in the following section using arguments similar to those used in the proof of Lemma \ref{LemmaMajArcFirst}. Let $Q_2 = n^{4/9}$, then we have the following minor arc estimate.

\begin{lemma}
\label{EstMFirstPruning}
One has
$$
\int_{\M(Q_1)\cap \m(Q_2)} |\cf_2^*(\a)|\text{d}\a \ll \cf(0) n^{\minus} n^{-\omega_e},
$$
with $\omega_e = 8.1\times10^{-5} $ and
$$
\cf^*_2(\a) = f_2^*(\a)F_3(\a) g_4(\a)\cdots g_{13}(\a)g_{17}(\a).
$$
\end{lemma}

The proof of this lemma follows a very similar approach to that of Lemma \ref{UpsNQ}, but the presence of $F_3$ unfortunately prevents us from using the same line of reasoning. The final pruning step is inspired by \cite[Lemma 2.5]{LiuZhao2021_RepresentationUnlikePowers}, where we will make crucial use of an estimate of the fourth-order moment of $F_3$ on major arcs \cite{Brudern1991}.

\begin{proposition}
\label{EstThirdPruning}
One has
$$
\int_{\M(Q_2)\cap\m(Q^{2/9})} |\cf_2^*(\a)|\text{d}\a \ll \cf(0)n^\minus n^{-\frac{1}{36}}.
$$
\end{proposition}

The remaining part will be routine and will be covered briefly, as most of the techniques used have already been discussed in Section \ref{MajArcsasymptotic}.

\subsection{Pruning stages}

As is now standard practice when implementing Ford's method, we take
$$L = \set{4,5,6,8,10} \et K = \set{7,9,11,12,13,17}.$$

\begin{lemma}
\label{EstGlGk}
One has
$$
\int_0^1 |F_3(\a)G_K(\a)|^2\text{d}\a \ll F_3(0)^2G_K(0)^2n^{-\phi_K},
$$
where $\phi_K = 0.796951$.
Similarly, we obtain
$$
\int_0^1 |G_L(\a)|^2\text{d}\a \ll G_L(0)^2n^{-\phi_L},
$$
where $\phi_L = 0.728613$.
Also, for any $1 \le Q \le n^\half$, we have the bound
$$
\int_{\M(Q)} |G_L(\a)|^2\text{d}\a \ll G_L(0)^2n^\minus Q^{\omega_L},
$$
where $\omega_L = 0.542744$
\end{lemma}

A slight difference between Section \ref{FirstMajArc} and this paragraph lies in the presence of the exponential sum $F_3$ within the integral. A line of reasoning similar to Ford’s method presented in Section \ref{SecMVE} also applies to estimates of mixed mean values,\textit{ i.e.} where $F_3$ appears in addition to a certain $g_k$ in the second moment order. Indeed, in line with Ford’s work \cite{Ford1996UnlikePowersII,KLX2022}, the same type of estimates as those presented in Section \ref{SecMVE} always allow us to recover appropriate bounds for mixed mean values. Indeed, for any choice of $K \sqcup L =\set{4,\dots,13}\cup\{17\}$, we have
\vspace{-2pt}
\begin{equation}
    \int_{\m(Q)}|F_3g_4(\a)\cdots g_{13}(\a)g_{17}(\a)|\text{d}\a \ll \left(\int_{0}^1|F_3G_K(\a)|^2\text{d}\a\right)^\half\left(\int_{0}^1|G_L(\a)|^2\text{d}\a\right)^\half.
\end{equation}

\leftcentersright[0]{Then, introduce}{$\displaystyle \hat\Kappa := \int_{0}^1|F_3(\a)G_K(\a)|^2\text{d}\a,$}{}

\leftcenters[0]{and}{$\displaystyle \hat{U}_{k,s}(P) := \int_0^1 |F_3(\a)g_k(\a)^s|^2 \text{d}\a.$}

\noindent In particular, \cite{Ford1996UnlikePowersII} exhibits a recursive algorithm implemented by \cite{KLX2022} to compute a certain $\omega(k,s)$ such that
\vspace{-1pt}
\begin{equation}
\label{Uhatks}
    \hat{U}_{k,s}(P) \ll (F_3(0)g_k(0)^s)^2n^{-\omega(k,s)},
\end{equation}

\noindent Again, performing a Hölder inequality of convexity with exponents $\{s_k\}_{k\in K}$ leads to
\vspace{-1pt}
$$
\hat\Kappa \ll \prod\limits_{k \in K_j} \hat{U}_{k,2s_k}(n^\frac{1}{k})^{\frac{1}{s_k}} \ll (F_3(0)G_K(0))^2 n^{-\phi_{\hat K}}
$$

\noindent where $\phi_{\hat K} = \sum\limits_{k \in K}\frac{\omega(k,s)}{s_k}$.

\begin{proof}{(\textit{\autoref{EstGlGk}})}
    Using the definition of $K \et L$, we define the set of Hölder exponents $s_k$ according to the values in \autoref{SetKL} where $s_{4} >0 \et s_{11} >0$ are chosen such that

\begin{equation*}
\frac{1}{4}+\frac{1}{6}+\frac{1}{s_{11}}+\frac{1}{6}+\frac{1}{6.75}+\frac{1}{10} = \frac{1}{s_4}+\frac{1}{4}+\frac{1}{5}+\frac{1}{7}+\frac{1}{8} = 1.    
\end{equation*}

\noindent We note that $s_{4} = 3.544303\dots$ and $s_{11} =  5.934066\dots$

\begin{table}[H]
\centering
\begin{minipage}{0.4\textwidth}
\centering
\begin{tabular}{|c|c|}
\hline
$k$ & $s_k$\\
\hline
   7 & 4\\
    9& 6\\
  11    & $s_{11}$\\
  12     & 6\\
  13 & 6.75 \\
  17     & 10\\
\hline
\end{tabular}
\\
\vspace{5pt}
Hölder's exponents for~$K$
\end{minipage}
\hspace{5pt}
\begin{minipage}{0.4\textwidth}
\centering
\begin{tabular}{|c|c|}
\hline
$k$ & $s_k$\\
\hline
4 & $s_4$\\
5 & 4\\
6 & 5\\
8 & 7\\
10 & 8 \\
\hline
\end{tabular}
\\
\vspace{5pt}
Hölder's exponents for~$L$
\end{minipage}
\caption{Hölder's exponents for $K \et L$}
\label{SetKL}
\end{table}

Therefore, using \autoref{TabAdmExpK7andK4} of admissible exponents, we arrive at the estimate :
$$
\int_0^1 |F_3(\a)G_K(\a)|^2\text{d}\a \ll (F_3(0)G_K(0))^2n^{-\phi_K}
$$
and $\phi_K = 0.796951$. Also, one has
$$
\int_0^1 |G_L(\a)|^2\text{d}\a \ll G_L(0)^2n^{-\phi_L}
$$
and $\phi_L = 0.728613$.
The first two results stated in Lemma \ref{EstGlGk} are therefore proven. The final result follows from the crucial observation in \eqref{DefAlphaK}. Using \autoref{Theorem 4.2}, combined with the Hölder exponents presented in the introduction to the proof, we get 

$$
\int_{\M(Q)} |G_L(\a)|^2\text{d}\a \ll G_L(0)^2n^\minus Q^{2(1-\phi_L)}.
$$

\noindent The value of $\phi_L$ obtained previously leads to the required conclusion, that is
$\omega_L = 2(1-\phi_L) = 0.542744$.
\end{proof}

Now using both Lemma \ref{EstGlGk} and \autoref{WeylEstimate} for $F_2$, we can prove Lemma \ref{EstMdelta}. \\ Let $\delta = \frac{1}{2}\left(\phi_K+\phi_L-\frac{3}{2}\right) = 12.782\times10^{-3}$ such that 

\begin{equation}
\label{Valdelta}
    \phi_K + \phi_L = \frac{3}{2}+2\delta.
\end{equation}

\noindent We apply to the integral over $\m(Q_1)$ a now well-known procedure. Whenever $\a \in \m(Q_1),$ \eqref{EstF2Minor} provides

$$
\int_{\m(Q_1)} |\cf(\a)|\text{d}\a \ll F_2(0)Q_1^{-\half + \varepsilon} \int_{\m(Q_1)}|F_3(\a)G_K(\a)G_L(\a)|\text{d}\a.
$$

\noindent Next, applying Cauchy's inequality and combining both results of \autoref{EstGlGk} yields

\begin{equation}
    \label{EqMDelta}
\int_{\m(Q_1)} |\cf(\a)|\text{d}\a \ll F_2(0)Q_1^{-\half+\varepsilon} F_3(0)G_K(0)n^{-\phi_K/2}G_L(0)n^{-\phi_L/2}.
\end{equation}

\noindent Therefore, as $Q_1 = n^{\half - \delta}$, recall \eqref{Valdelta} and use \eqref{EqMDelta} to conclude that
$$
\int_{\m(Q_1)} |\cf(\a)|\text{d}\a \ll \cf(0)n^{-1+\delta/2-\delta+\varepsilon} \ll \cf(0) n^\minus n^{-\delta/2 +\varepsilon}.
$$

Before continuing with our proof, let us briefly note that the major arc we have just reached has exponent $1/2-\delta \simeq 0.487218$. This can be compared with the minor arcs $\m(n^\mu)$ with $\mu \simeq 0.476675$ used by Liu and Zhao in \cite[Lemma 2.3] {LiuZhao2021_RepresentationUnlikePowers}. This is the best value obtainable with the current data in the various tables of admissible exponents. Now that we have reached a non-trivial major arc, we can consider replacing $F_2$ with its associated Weyl sum, as in Lemma \ref{LemmaReplacef2}.

\begin{lemma}
\label{EstReplaceF2}
One has
$$
\int_{\M(Q_1)} |\Delta_2(\a)F_3(\a)G_K(\a)G_L(\a)|\text{d}\a \ll \cf(0) n^\minus n^{-\delta+\varepsilon}.
$$
\end{lemma}

\begin{proof}{}

The last remark we have just made allows us to reproduce, to a large extent, the proof of \autoref{LemmaReplacef2}. Let us begin by bounding the term $\Delta_2(\a)$ by $Q_1^{\half + \varepsilon}$. Next, using Cauchy’s inequality and completing the integral over $\M(Q_1)$ with an integral over $[0,1)$, we apply \autoref{EstGlGk} to obtain

$$
\int_{\M(Q_1)} |\Delta_2(\a)F_3(\a)G_K(\a)G_L(\a)|\text{d}\a \ll Q_1^{\half + \varepsilon}F_3(0)G_K(0)G_L(0) n^{-\phi_K/2-\phi_L/2}.
$$

\noindent Finally, we substitute $Q_1$ with its value and introduce $F_2(0)$ into the previous equation to get an upper bound that is slightly better than that given in Lemma \ref{EstReplaceF2}. This concludes the proof of the desired result.
\end{proof}

We can now briefly summarise what we have established so far. According to the results of Lemmata \ref{EstMdelta} and \ref{EstReplaceF2}, we have
$$
H(n;\eta) = \int_{\M(Q_1)} \cf_2^*(\a) e(-n\a)\text{d}\a + \bigoh{\cf(0)n^\minus n^{-\delta}},
$$
where
$$
\cf_2^* := F_2^*F_3G_KG_L \et \delta >0.
$$

\noindent We now have the key ingredients to start the pruning process.

\begin{proof}{(\textit{\autoref{EstMFirstPruning}})}
The framework for proving this lemma is based on a dyadic dissection of $\M(Q_1) \cap \m(Q_2)$, where $Q_2 = n^{4/9}$, as defined in Lemma \ref{EstMFirstPruning}. We then have

$$
\M(Q_1) \cap \m(Q_2) = \bigcup\limits_{j=0}^J \cn(2^{-j}Q_1),
$$

\noindent where $J = \floor{\frac{\log(Q_1/Q_2)}{\log 2}} \ll \log n$. Having established this, we must restrict our proof to an integral over $\cn(Q)$ for $Q$ lying between $Q_2$ and $Q_1$.
First, we use \eqref{EstF2*} and the Cauchy inequality to obtain

$$
\int_{\cn(Q)} |\cf_2^*(\a)|\text{d}\a \ll F_2(0)Q^{-\half+\varepsilon}\left(\int_0^1 |F_3(\a)G_K(\a)|^2\text{d}\a \right)^\half \left(\int_{\M(Q)} |G_K(\a)|^2\text{d}\a\right)^\half.
$$

\noindent Consequently, Lemma \ref{EstGlGk} yields the estimate

$$
\int_{\cn(Q)} |\cf_2^*(\a)|\text{d}\a \ll \cf(0) n^{-\frac{\phi_L}{2}}n^{-\half}Q^{\frac{\omega_K}{2}-\half+\varepsilon}.
$$

\noindent But since we are in the range where $Q \ge n^{4/9}$, we have
\begin{align*}
    \int_{\cn(Q)} |\cf_2^*(\a)|\text{d}\a &\ll \cf(0) n^{-\frac{\phi_L}{2}-\half+\varepsilon}n^{\frac{4}{9}(-\half+\frac{\omega_K}{2})} \\
    &\ll \cf(0) n^{-1}n^{-\omega_e+\varepsilon} ,
\end{align*}

\noindent where $\omega_e = -1+\left(\frac{\phi_L}{2}+\frac{1}{2}-\frac{4}{9}\left(-\frac{1}{2}+\frac{\omega_K}{2}\right)\right) \simeq 8.1\times10^{-5}.$
\end{proof}

Now that we have reached $\M(n^{\frac{4}{9}})$, we need to adjust our approach to the pruning process. Essentially, we have made the most of the continuous gain in the savings over major arcs based on smooth estimates of the mean value given by \autoref{EstGlGk}. As in the present case, the sum $F_3$ precludes such reasoning; we can compare Lemma \ref{EstGlGk} with Lemma \ref{UpsNQ}, and we follow the ideas in Section $5$ of \cite{LiuZhao2021_RepresentationUnlikePowers} to prune major arcs.

\noindent In any case, since we are now going to begin decoupling $F_3 \et G_K$, we still need to determine sufficiently good estimates of the mean value of $G_KG_L$ on various sets. We therefore introduce a new partition of $\set{4,\dots,13,17}$. Let us write
$$
\widetilde{K} = \set{4,8,9,10,11,12,13} 
\et \widetilde{L} = \set{5,6,7,17}.
$$

\begin{lemma}
\label{EstGLtGKt}
Let $1 \le Q \le n^\half$. One has
$$
\int_{\M(Q)} |G_{\widetilde{L}}(\a)|^4\text{d}\a \ll G_{\widetilde{L}}(0)^4n^\minus Q^{\omega_{\widetilde{L}}},
$$
where $\omega_{\widetilde{L}} = 0.309666$, and
$$
\int_{\M(Q)} |G_{\widetilde{K}}(\a)|^2\text{d}\a \ll G_{\widetilde{L}}(0)^2n^\minus Q^{\omega_{\widetilde{K}}},
$$
with $\omega_{\widetilde{K}} = 
0.591714$.
\end{lemma}

\begin{proof}{}
Using $\widetilde{K} \et \widetilde{L}$ as in \autoref{EstGLtGKt}, an exhaustive application of Ford's algorithm, with Hölder exponents and tabulated values depicted in \autoref{TabAdmExpGKtGLt}, yields 
\vspace{-1pt}
\begin{align}
\label{EqGLt}
    \int_{\M(Q)} |G_{\widetilde{L}}(\a)|^4\text{d}\a 
    &\ll \left(\int_{\M(Q)}|g_5(\a)^4g_7(\a)g_{17}(\a)^3|^2\text{d}\a \right)^\half\left(\int_{\M(Q)}|g_6(\a)^4g_7(\a)^3g_{17}(\a)|^2\text{d}\a \right)^\half \nonumber \\
    &\ll V_{5,2s}(Q)^\frac{4}{2s}  V_{7,16}(Q)^\frac{1}{16} V_{17,44}(Q)^\frac{3}{44}  V_{6,14}(Q)^\frac{4}{14} V_{7,2t}(Q)^\frac{3}{2t} V_{17,40}(Q)^\frac{1}{40} \nonumber \\
    &\ll G_{\widetilde{L}}(0)^4 n^\minus Q^{2(1-\phi_{\widetilde{L}})}
\end{align}

\leftcenters[0]{with}{$s  = \frac{4}{1-(\frac{1}{8}+\frac{3}{22})} \simeq 5.4153 \et t = \frac{3}{1-(\frac{4}{7}+\frac{1}{20})} \simeq 7.9245.$}

\noindent Then, as announced, we achieve the exponents 

$$\phi_{\widetilde{L}} = 0.845166 \et  \omega_{\widetilde{L}} = 2(1-\phi_{\widetilde{L}}) = 0.309666.$$

\noindent Regarding the second moment estimate for $G_{\widetilde{K}}$, the same method with \autoref{TabAdmExpGKtGLt} yields
\vspace{-1pt}
\begin{align}
\label{EqGKt}
    \int_{\M(Q)}|G_{\widetilde{K}}(\a)|^2\text{d}\a 
    &\ll \left(\int_{\M(Q)}|g_8(\a)g_9(\a)^2g_{10}(\a)g_{11}(\a)g_{12}(\a)^2g_{13}(\a)^2|^2\text{d}\a \right)^\half \\ &\quad \quad \quad \left(\int_{\M(Q)}|g_4(\a)^2g_8(\a)g_{10}(\a)g_{11}(\a)|^2\text{d}\a \right)^\half \nonumber \\
    &\ll V_{8,14}(Q)^\frac{1}{14}  V_{9,16}(Q)^\frac{2}{16} V_{10,2s}(Q)^\frac{1}{2s} V_{11,18}(Q)^\frac{1}{18} V_{12,20}(Q)^\frac{2}{20} V_{13,22}(Q)^\frac{2}{22} \\ & \qquad \qquad V_{4,2t}(Q)^\frac{2}{2t}V_{8,12}(Q)^\frac{1}{12} V_{10,16}(Q)^\frac{1}{16} V_{11,18}(Q)^\frac{1}{18} \nonumber \\
    &\ll G_{\widetilde{K}}(0)^2n^\minus Q^{\omega_{\widetilde{K}}},
\end{align}

\leftcenters[0]{where}{$s = \frac{1}{1-(\frac{1}{7}+\frac{2}{8}+\frac{1}{9}+\frac{2}{10}+\frac{2}{11})} \simeq 8.755527; \quad  t = \frac{2}{1-(\frac{1}{6}+\frac{1}{8}+\frac{1}{9})} \simeq 3.348837 \et \omega_{\widetilde{K}} = 0.591714.$} 
\end{proof}

Next, we recall a lemma due to Br\"{u}dern.

\begin{lemma}{\cite[Theorem 2]{Brudern1991}}
\label{MVF3}
For any $1 \le Q \le n^{1/3}$, one has
$$
\int_{\M(Q)} |F_3(\a)|^4\text{d}\a \ll n^\varepsilon(n^{1/3}+Q^2n^{-1/3}+Q^{7/2}n^\minus).
$$
\end{lemma}

Note that for $n^{1/3} \le Q \le n^{4/9}$, the dominant term in the previous lemma is the middle one. And for $Q \le n^{1/3}$, the dominant term is the one on the far left, which provides a quasi-optimal upper bound for this fourth moment of $F_3$.
All this enables us to secure a pruning up to $Q_3 = n^{1/3}$.

\begin{lemma}
\label{EstScndPruning}
Since $Q_2 \le n^{4/9}$, we have
$$
\int_{\M(Q_2)\cap\m(Q_3)} |\cf_2^*(\a)|\text{d}\a \ll \cf(0)n^\minus n^{-\delta}
$$
for some positive $\delta$.
\end{lemma}

\begin{proof}{}
Let $n^{1/3} \le Q \le n^{4/9}$.
As stated in Lemma \ref{MVF3}, in this range of values for $Q$, we have
$$
\int_{\M(Q)} |F_3(\a)|^4\text{d}\a \ll n^\varepsilon Q^2n^{-1/3}.
$$
We may use a dyadic dissection of $\M(n^{4/9})\cap \m(n^{1/3})$ to conclude here. Using Hölder inequality and \eqref{EstF2*}, observe that
\vspace{-2pt}
$$
\int_{\cn(Q)} |\cf_2^*(\a)|\text{d}\a \ll F_2(0)Q^{-\half + \varepsilon} \left(\int_{\M(Q)} |F_3(\a)|^4\text{d}\a\right)^\frac{1}{4} \left(\int_{\M(Q)}|G_{\widetilde{L}}(\a)|^4\text{d}\a\right)^{\frac{1}{4}} \left(\int_{\M(Q)}|G_{\widetilde{K}}(\a)|^2\text{d}\a\right)^{\frac{1}{2}}.
$$

\noindent Then, using both results of Lemma \ref{EstGLtGKt} yields
\vspace{-1pt}
\begin{align}
\label{intermediate29}
    \int_{\cn(Q)} |\cf_2^*(\a)|\text{d}\a 
    &\ll F_2(0)Q^{-\half+\varepsilon}\left(n^\varepsilon Q^2 n^{-1/3}\right)^{\frac{1}{4}} \left(G_{\widetilde{L}}(0)^4 n^\minus Q^{\omega_{\widetilde{L}}}\right)^{\frac{1}{4}}\left(G_{\widetilde{K}}(0)^2 n^\minus Q^{\omega_{\widetilde{K}}}\right)^\half  \nonumber\\
    &\ll F_2(0)F_3(0)n^{-\frac{1}{3}}G_{\widetilde{K}}(0)G_{\widetilde{L}}(0)n^{-\frac{5}{6}}Q^{\frac{\omega_{\widetilde{L}}}{4}+\frac{\omega_{\widetilde{K}}}{2}+\varepsilon} \nonumber\\
    &\ll \cf(0) n^\minus Q^{\frac{\omega_{\widetilde{L}}}{4}+\frac{\omega_{\widetilde{K}}}{2}+\varepsilon}n^{-\frac{1}{6}}.
\end{align}

Finally, as $n^\frac{1}{3}\le Q \le n^\frac{4}{9}$, note that $\frac{4}{9}\left(\frac{\omega_{\widetilde{L}}}{4}+\frac{\omega_{\widetilde{K}}}{2}\right) = 0.165899\dots < 1/6$. Writing $\delta = 1/6-\frac{4}{9}\left(\frac{\omega_{\widetilde{L}}}{4}+\frac{\omega_{\widetilde{K}}}{2}\right) $, we observe that the rightmost part of \eqref{intermediate29} is

$$
\int_{\cn(Q)} |\cf_2^*(\a)|\text{d}\a \ll \cf(0) n^\minus n^{-\delta + \varepsilon}.
$$

\noindent Then, proceeding to a dyadic dissection of $\M(n^{4/9}) \backslash \M(n^{1/3})$, one has

$$
\int_{\M(n^{4/9}) \backslash \M(n^{1/3})} |\cf_2^*(\a)|\text{d}\a \ll \cf(0) n^\minus n^{-\delta + \varepsilon}.
$$

\end{proof}

Once we have reached the main arc $\M(n^{1/3})$, of a reasonable size, we can carry out further pruning to remain within the validity range of \eqref{EstimateDelta3}, which is precisely the aim of Lemma \ref{EstThirdPruning}.

\begin{proof}{\textit{(Lemma \ref{EstThirdPruning})}}
    Let $n^{2/9} \le Q \le n^{1/3}$. Then, note that
    
    $$
    \int_{\M(Q)}|F_3(\a)|^4\text{d}\a \ll F_3(0)^4n^{-1+\varepsilon}.
    $$
    
    \noindent As this bound is essentially optimal (up to a $n^\varepsilon$ term), we imitate the proof of Lemma \ref{EstScndPruning} to achieve 
    \vspace{-1pt}
    $$
    \int_{\cn(Q)} |\cf_2^*(\a)|\text{d}\a \ll F_2(0)Q^{-\half+\varepsilon}\left(F_3(0)^4n^\minus n^\varepsilon\right)^\frac{1}{4}\left(G_{\widetilde{L}}(0)^4 n^\minus Q^{\omega_{\widetilde{L}}}\right)^{\frac{1}{4}}\left(G_{\widetilde{K}}(0)^2n^\minus Q^{\omega_{\widetilde{K}}}\right)^\half.
    $$
    
    \noindent Thus, gathering all interesting terms altogether, and introducing $\delta_2 = 0.126725 < \half - (\frac{\omega_{\widetilde{L}}}{4}+\frac{\omega_{\widetilde{K}}}{2})$, we reach the bound
    \vspace{-2pt}
    \begin{align*}
        \int_{\cn(Q)}|\cf_2^*(\a)|\text{d}\a 
        &\ll \cf(0) n^{\minus+\varepsilon} Q^{-\delta_2+\varepsilon}.
    \end{align*}

    \noindent Since $\delta_2 > \frac{1}{8}$, performing a dyadic dissection of $\M(n^{1/3})\cap \m(n^{2/9})$ leads to
    $$
    \int_{\M(n^{1/3})\cap \m(n^{2/9})}|\cf_2^*(\a)|\text{d}\a \ll \cf(0)n^\minus n^{-1/36}.
    $$
\end{proof}

We may conclude this section with the following remark, which summarises the essence of the previous pruning. When $\a \in \M(n^{\frac{1}{3}-\delta})$ and for $k = 3$, note that we reach a non-trivial estimate of $\Delta_3$. Indeed, a partial summation formula applied to \eqref{Ek} allows us to generalise the result for $k = 3$. Indeed, for any $1 \le Q \le n^{\frac{1}{3}-\delta}$ and $\a \in \M(Q)$, \cite[Theorem 4.1]{vaughan1981hardy} gives

\begin{equation}
    \label{EstimateDelta3}
\Delta_3(\a) \ll q^{\half+\varepsilon} \ll Q^{\half + \varepsilon}.
\end{equation}

\noindent And when $\a \in \frakN(Q)$, again \cite[Lemma 2.3 and 2.8]{vaughan1981hardy} yield

\begin{equation}
\label{EstF3*}
    F_3^*(\a) \ll F_3(0)\frac{q^{-\frac{1}{3}+\varepsilon}}{1+n|\b|} \ll F_3(0)Q^{-\frac{1}{3}+\varepsilon}.
\end{equation}

\begin{lemma}
\label{MVGKGL}

Let $1 \le Q \le n^\half$, one has

$$
\int_{\M(Q)} |G_K(\a)G_L(\a)|\text{d}\a \ll G_K(0)G_L(0)n^\minus Q^{\omega_{KL}},
$$

\noindent where $\omega_{KL} = 0.748615$. Therefore, we obtain

$$
\int_{\M(n^{2/9})} |F_2^*(\a)\Delta_3(\a)G_K(\a)G_L(\a)|\text{d}\a \ll \cf(0)n^\minus n^{-\frac{1}{6}}.
$$

\end{lemma}

\begin{proof}{}

    Now, to replace $F_3$ by $F_3^*$, we have to study
    $$
    \int_{\M(n^{2/9})}|F_2^*(\a)\Delta_3(\a)G_K(\a)G_L(\a)| \text{d}\a.
    $$
    Using \eqref{EstF2*} and \eqref{EstimateDelta3}, we first bound this integral by
    
    \begin{equation}
        \label{Process}
        F_2(0)F_3(0)n^{-1/3+\varepsilon}\int_{\M(n^{2/9})} |G_K(\a)G_L(\a)|\text{d}\a.
    \end{equation}
    
    All that remains is to obtain an estimate of the integral. As this will be of great interest in the next pruning step, applying the same approach as in \autoref{EstGLtGKt}, we claim the following general estimate for all $1 \le Q \le n^\half$. Using \autoref{TabAdmExpGKGL} for the set of exponents and the numerical data, we get
    \vspace{-2pt}
    \begin{align}
    \label{EqGKGL}
        \int_{\M(Q)}|G_K(\a)G_L(\a)|\text{d}\a &\ll G_K(0)G_L(0)n^\minus Q^{\omega_{KL}},
    \end{align}
    
    \noindent where $\omega_{KL} = 0.748615$ is defined in \autoref{MVGKGL}.
    Consequently, using this result in \eqref{Process}, we obtain
    \vspace{-3pt}
    \begin{align*}
        \int_{\M(n^{2/9})}|F_2^*(\a)\Delta_3(\a)G_K(\a)G_L(\a)|\text{d}\a 
        &\ll \cf(0)n^\minus n^{-\frac{1}{3}} n^{\frac{2}{9}\cdot \omega_{KL}+\varepsilon} \\
        &\ll \cf(0)n^\minus n^{-\frac{1}{6}}.
    \end{align*}
\end{proof}

We are now ready to deliver the decisive step towards achieving a more manageable major arc. As mentioned in Lemma \ref{MVGKGL}, we can substitute $F_3$ with $F_3^*$ and define

$$
\cf_3^* := F_2^*F_3^*G_KG_L.
$$

\noindent Let $1 \le Q \le n^{2/9}$. Observe that over $\cn(Q)$, \eqref{EstF2*} and \eqref{EstF3*} ensure that

$$
F_2^*(\a)F_3^*(\a) \ll F_2(0)F_3(0) Q^{-5/6+\varepsilon}.
$$

\noindent Therefore, this estimate combined with \autoref{MVGKGL} yields
\vspace{-1pt}
\begin{align*}
\int_{\cn(Q)} |\cf_3^*(\a)|\text{d}\a
&\ll F_2(0)F_3(0)Q^{-5/6+\varepsilon}\int_{\M(Q)}|G_K(\a)G_L(\a)|\text{d}\a \\
&\ll \cf(0) n^\minus Q^{-\frac{5}{6}+\omega_{KL}+\varepsilon}.
\end{align*}

\noindent Hence
\vspace{-2pt}
\begin{equation}
    \label{LastPruning}
    \int_{\cn(Q)} |\cf_3^*(\a)|\text{d}\a \ll \cf(0)n^\minus Q^{-\frac{1}{12}}.
\end{equation}

Finally, recall that $\Kfrak$ is defined as the union of the intervals

$$
\Kfrak(q,a) = \set{\a \in [0,1) : |\a - \frac{a}{q}| \le \frac{L^\frac{1}{15}}{n}},
$$

\noindent with $0 \le a \le q \le L^\frac{1}{15} \et (a,q)=1$. As was done in \autoref{Lemmashrink}, we can reach an integral over $\Kfrak$ by using~\eqref{LastPruning}
\vspace{-4pt}
\begin{align}  
\label{EndofPruning}
\int_{\M(n^{2/9})\backslash \Kfrak} |\cf_{3}^*(\a)|\text{d}\a
&\ll \int_{\M(n^{2/9})\cap \m(L^{1/15})} |\cf_{3}^*(\a)|\text{d}\a \nonumber \\
& \ll \sum\limits_{2^{-j}n^{2/9}\ge L^{1/15}} \int_{\cn(2^{-j}n^{2/9})}|\cf_3^*(\a)|\text{d}\a \nonumber \\
&\ll \cf(0) n^\minus \sum\limits_{2^{-j}n^{2/9} \ge L^{1/15}} \left(2^{-j}n^{2/9}\right)^{-1/12}\nonumber \\
&\ll \cf(0)n^\minus L^{-1/12\cdot1/15}.
\end{align}

\noindent Let $\rho = 1/180$, then all the previous Lemmata and \eqref{EndofPruning} fit together to form

\begin{equation}
\label{FinalBlow}
H(n) = \int_{\Kfrak} \cf_{3}^*(\a)e(-n\a)\text{d}\a +\bigoh{\cf(0)n^\minus L^{-\rho}}.
\end{equation}

\noindent As \eqref{FinalBlow} summarises, we are left with the following integral to study

$$
\int_{\Kfrak} \cf_3^*(\a)e(-n\a)\text{d}\a.
$$

\noindent The final step before moving on to the asymptotic case is to replace the smooth Weyl sums with their smooth ‘cousin’, as suggested by Lemma \ref{SmoothWeilCousin}. In line with the proof of Lemma \ref{ReplaceSmoothCousin}, we can replace each $g_k$ one by one to finally reach

$$
\int_{\Kfrak} \cf_3^*(\a)e(-n\a)\text{d}\a  = \int_{\Kfrak} \cf^*(\a)e(-n\a) \text{d}\a + \bigoh{\cf(0)L^{-\frac{1}{20}}n^\minus},$$

\noindent Once this result is secured, we introduce the singular series in the form
\vspace{-3pt}
\begin{align*}
A(n,q) &= c^{11}\frac{1}{q^{13}}\sum\limits_{(a,q)=1}S_2\cdots S_{13}S_{17}(q,a)e(-an/q) \\
\et \S(n,Y) &= \sum\limits_{q \le Y} A(n,q).
\end{align*}

Let $\I(n)$ denote the singular integral defined as 
\vspace{-2pt}
\begin{align*}
\I(n,Y) &= \int_{-Y/n}^{Y/n} \omega_2(x)w_3(x)\prod\limits_{k=4}^{13} \nu_k(x)\nu_{17}(x)dx, \\
\et \I(n) &= \I(n,1/2).
\end{align*}

\noindent Therefore, a straightforward application of the methods used in Sections \ref{SecSingularSeries} and \ref{SecSingularIntegral} provide

\begin{equation*}
    H(n;\eta) = \S(n)\I(n) + \bigoh{n^{\theta-1}L^{-\tau}}.
\end{equation*}

\noindent Finally, it is not difficult to show that $\S(n)\I(n) \gg n^{\theta-1}$, and the conclusion follows.

\section{Improving a previous result of Ford}
\label{StartThree}

We finish by sketching a proof of the following result, which improves Theorem 2 of \cite{Ford1996UnlikePowersII}.

\ThUnlikePowersThreeresta*

\noindent Building on what we have covered in the previous sections, let us take 
\begin{align*}
    K_4 = \{&4,6,7,9,15,17,20,21,22,24,28,29,32,34,35,37,39,\\
    &40,41,43,45,47,51,54,58,61,63,64,65,66,69,70,71\}, \\
    \et K_5 = \{&5,8,10,11,12,13,14,16,18,19,23,25,26,24,30,31,33,36,38,\\
    &42,44,46,48,49,50,52,53,55,56,57,59,60,62,67,68,72,73\}.
\end{align*}

\noindent with sets of exponents presented in \autoref{SetK4again} and \autoref{SetK5} below.

\begin{table}[H]
\centering

\begin{minipage}{0.23\textwidth}
\centering
\begin{tabular}{|c|c|}
\hline
$k$ & $s_k$\\
\hline
4 & $s_4$\\
6 & 8\\
7 & 10\\
9 & 13\\
15 & 24 \\
17 & 28 \\
20 & 33 \\
21 & 34 \\
22 & 36 \\
\hline
\end{tabular}
\end{minipage}
\hspace{-25pt}
\begin{minipage}{0.23\textwidth}
\centering
\begin{tabular}{|c|c|}
\hline
$k$ & $s_k$\\
\hline
24 & 39 \\
28 & 46 \\
29 & 48 \\
32 & 53 \\
34 & 56 \\
35 & 57 \\
37 & 61 \\
39 & 64 \\
\hline
\end{tabular}
\end{minipage}
\hspace{-25pt}
\begin{minipage}{0.23\textwidth}
\centering
\begin{tabular}{|c|c|}
\hline
$k$ & $s_k$\\
\hline
40 & 66 \\
41 & 67 \\
43 & 71 \\
45 & 73 \\
47 & 77 \\
51 & 84 \\
54 & 89 \\
58 & 95 \\
\hline
\end{tabular}
\end{minipage}
\hspace{-25pt}
\begin{minipage}{0.23\textwidth}
\centering
\begin{tabular}{|c|c|}
\hline
$k$ & $s_k$\\
\hline
61 & 100 \\
63 & 103 \\
64 & 105 \\
65 & 107 \\
66 & 108 \\
69 & 113 \\
70 & 115 \\
71 & 116 \\
\hline
\end{tabular}
\end{minipage}

\caption{Hölder's exponents for $K_4$}
\label{SetK4again}
\end{table}

\begin{table}[H]
\centering

\begin{minipage}{0.23\textwidth}
\centering
\begin{tabular}{|c|c|}
\hline
$k$ & $s_k$\\
\hline
5 & 7 \\
8 & 12 \\
10 & 15 \\
11 & 17 \\
12 & 18 \\
13 & 20 \\
14 & 23 \\
16 & $s_{16}$ \\
18 & 29 \\
19 & 31 \\
\hline
\end{tabular}
\end{minipage}
\hspace{-25pt}
\begin{minipage}{0.23\textwidth}
\centering
\begin{tabular}{|c|c|}
\hline
$k$ & $s_k$\\
\hline
23 & 37 \\
25 & 41 \\
26 & 42 \\
27 & 44 \\
30 & 48 \\
31 & 51 \\
33 & 54 \\
36 & 59 \\
38 & 62 \\
\hline
\end{tabular}
\end{minipage}
\hspace{-25pt}
\begin{minipage}{0.23\textwidth}
\centering
\begin{tabular}{|c|c|}
\hline
$k$ & $s_k$\\
\hline
42 & 68 \\
44 & 72 \\
46 & 75 \\
48 & 78 \\
49 & 80 \\
50 & 81 \\
52 & 85 \\
53 & 86 \\
55 & 89 \\
\hline
\end{tabular}
\end{minipage}
\hspace{-25pt}
\begin{minipage}{0.23\textwidth}
\centering
\begin{tabular}{|c|c|}
\hline
$k$ & $s_k$\\
\hline
56 & 91 \\
57 & 93 \\
59 & 96 \\
60 & 98 \\
62 & 101 \\
67 & 109 \\
68 & 110 \\
72 & 117 \\
73 & 119 \\
\hline
\end{tabular}
\end{minipage}
\caption{Hölder's exponents for $K_5$}
\label{SetK5}
\end{table}

In which $s_4 = 4.995380\dots \et s_{16} = 24.479295\dots$ are precisely defined to get Hölder sets as in Tables \ref{SetK4} and \ref{SetK10}. Thus, Tables \ref{SetK4again} and \ref{SetK5}, combined with Tables \ref{TabAdmExpR31K4} and \ref{TabAdmExpR31K5}, yield the following result.

\begin{lem}
\label{EqR31}
    One has
    $$
    \int_{\cn(Q)} |G_4(\a)|^2\text{d}\a \ll G_4(0)^2n^\minus Q^{\omega_4},
    $$
    where $\omega_4 = 0.161975$, and
    $$
    \int_{\cn(Q)} |G_5(\a)|^2\text{d}\a \ll G_5(0)^2n^\minus Q^{\omega_5},
    $$
    where $\omega_5 = 0.170086$.
\end{lem}

Let us introduce $\delta = 3.18 \times 10^{-4} < 1/12 - \frac{\omega_4 + \omega_5}{4} $. We can therefore apply the same reasoning as in sections \ref{FirstMajArc} and \ref{Pruning}, replacing $f_2$ with $f_3$ and using Weyl’s inequality for $f_3$. This considerably increases the number of variables, since we can no longer save a power of $n^{1/4}$, but only a power of $n^{1/12}$, as for $ Q \le n^{1/3}$, it stems from \autoref{WeylEstimate} that
\vspace{-1pt}
$$
\sup\limits_{\a \in \m(Q)}|F_3(\a)|  \ll F_3(0)Q^{-\frac{1}{4}+\varepsilon} \ll F_3(0)n^{-\frac{1}{12}+\varepsilon}.
$$

\noindent It therefore follows that
$$
\int_{\m(n^{1/3})}|F_3(\a)G_4(\a)G_5(\a)|\text{d}\a \ll \sup\limits_{\a \in \m(n^{1/3})}|F_3(\a)|\left(\int_0^1 |G_4(\a)|^2\text{d}\a\right)^\half\left(\int_0^1 |G_5(\a)|^2\text{d}\a\right)^\half.
$$

\noindent It will also be interesting to use \eqref{EstimateDelta3} and \eqref{EstF3*} to replace $F_3$, as was done previously for $F_2$ in Lemma \ref{LemmaReplacef2}. 
\begin{lem}
    One has
    $$
    \int_{\m(n^{1/3})}|F_3(\a)G_4(\a)G_5(\a)|\text{d}\a \ll \cf(0)n^\minus  n^{-1/12+\omega_4/4+\omega_5/4+\varepsilon}.
    $$
    And therefore, for any $1 \le Q \le n^{1/3}$, one has
    $$
    \int_{\cn(Q)}|F_3^*(\a)G_4(\a)G_5(\a)|\text{d}\a \ll \cf(0)n^\minus Q^{-\delta+\varepsilon}.
    $$
\end{lem}

\noindent So we conclude in a similar fashion that 

$$
R(3,1) \le 71.
$$

\noindent However, it should be noted that we only reach the following

$$
\sum\limits_{k=3}^{73} \frac{1}{k} = 3.3745\dots
$$

\noindent which is very high compared with the value of $2.48$ previously claimed in Theorem \ref{ThEvenUnlikePowers}, or even with $2$, which is the theoretical limit of the circle method in this context.

\newpage
\appendix

\section{Admissible exponents}
\label{AppAdmExp}

\begin{table}[H]
\centering

\begin{minipage}{0.23\textwidth}
\centering
\begin{tabular}{|c|c|c|}
\hline
$k$ & $s_k$ & $\Delta_k$\\
\hline
4 & 3.65 & 0.761416\\
4 & 3.7 & 0.734081\\
6 & 6 & 1.231563\\
8 & 7 & 2.310599\\
22 & 22 & 5.901869\\
30 & 30 & 8.141245\\
34 & 35 & 8.830801\\
40 & 40 & 10.933884\\
50 & 50 & 13.723309\\
    52 & 52 & 14.280971\\
    60 & 60 & 16.511131\\
    62 & 62 & 17.068572\\
\hline
\end{tabular}
\end{minipage}
\hspace{15pt}
\begin{minipage}{0.23\textwidth}
\centering
\begin{tabular}{|c|c|c|}
\hline
$k$ & $s_k$ & $\Delta_k$\\
\hline

    66 & 67 & 17.751366\\
    68 & 70 & 17.884345\\
    70 & 72 & 18.441233\\
    74 & 76 & 19.555019\\
    88 & 90 & 23.453333\\
    90 & 92 & 24.010241\\
    92 & 93 & 24.994426 \\
    94 & 95 & 25.551517 \\
    96 & 97 & 26.108601\\
    100 & 102 & 26.794797\\
    106 & 107 & 28.893936\\
\hline
\end{tabular}
\end{minipage}
\hspace{15pt}
\begin{minipage}{0.23\textwidth}
\centering
\begin{tabular}{|c|c|c|}
\hline
$k$ & $s_k$ & $\Delta_k$\\
\hline
    
    108 & 110 & 29.022457\\
    110 & 111 & 30.008035\\
    118 & 120 & 31.807044\\    
    120 & 122 & 32.363963\\
    130 & 132 & 35.148564 \\
    132 & 134 & 35.705485\\
    134 & 136 & 36.262406\\
    144 & 146 & 39.047017\\
    152 & 153 & 41.705315\\
    154 & 156 & 41.831633\\
    156 & 158 & 42.388557\\
\hline
\end{tabular}
\end{minipage}
    \caption{Admissible exponents for \eqref{Kappa4} used with the set of \autoref{SetK4}}
    \label{TabAdmExpK4}
\end{table}

\begin{table}[H]
    \centering

\begin{minipage}{0.23\textwidth}
\centering
\begin{tabular}{|c|c|c|}
\hline
$k$ & $s_k$ & $\Delta_k$\\
\hline
10 & 9 & 2.894571\\
10 & 10 & 2.437568\\
12 & 12 & 3.017477\\
14 & 14 & 3.589278\\
16 & 16 & 4.157379\\
18 & 18 & 4.723122\\
20 & 20 & 5.285739\\
24 & 24 & 6.462444\\
26 & 26 & 7.022457\\
28 & 29 & 7.155856\\
32 & 32 & 8.700175\\
32 & 33 & 8.272650\\
36 & 37 & 9.388821\\
36 & 38 & 8.975115\\
38 & 39 & 9.946732\\
42 & 43 & 11.062282\\
\hline
\end{tabular}
\end{minipage}
\hspace{15pt}
\begin{minipage}{0.23\textwidth}
\centering
\begin{tabular}{|c|c|c|}
\hline
$k$ & $s_k$ & $\Delta_k$\\
\hline
44 & 45 & 11.619945\\
44 & 46 & 11.202112\\
46 & 47 & 12.177547\\
46 & 48 & 11.758913\\
48 & 49 & 12.735094\\
54 & 55 & 14.407467\\
56 & 58 & 14.543097\\
58 & 60 & 15.099959\\
64 & 66 & 16.770578\\
72 & 74 & 18.998125\\
76 & 78 & 20.111915\\
78 & 80 & 20.668814\\
80 & 83 & 20.806341\\
82 & 84 & 21.782617\\
84 & 86 & 22.339521\\
86 & 88 & 22.896426\\
98 & 101 & 25.815503\\
\hline
\end{tabular}
\end{minipage}
\hspace{15pt}
\begin{minipage}{0.23\textwidth}
\centering
\begin{tabular}{|c|c|c|}
\hline
$k$ & $s_k$ & $\Delta_k$\\
\hline

102 & 104 & 27.351711\\
104 & 106 & 27.908626\\
112 & 115 & 29.712243\\
114 & 117 & 30.268956\\
116 & 120 & 30.405782\\
122 & 125 & 32.495880\\
124 & 127 & 33.052628\\
126 & 129 & 33.609381\\
128 & 131 & 34.166140\\
136 & 139 & 36.393226\\
138 & 143 & 36.109014\\
140 & 145 & 36.665370\\
142 & 147 & 37.221741\\
146 & 150 & 38.754050\\
148 & 151 & 39.733985\\
150 & 155 & 39.447376\\
158 & 162 & 42.093964\\
\hline
\end{tabular}
\end{minipage}

    \caption{Admissible exponents for \eqref{Kappa10} used with the set of \autoref{SetK10}}
    \label{TabAdmExpK10}
\end{table}


\begin{table}[H]
    \centering
    \begin{minipage}{0.4\textwidth}
\centering
\begin{tabular}{|c|c|c|}
\hline
$k$ & $s_k$&$\omega_{k,s_k}$\\
\hline
    7 & 4 & 2.979764\\
    9 & 6 & 3.442305\\
    11 & 5 & 2.569268\\
    11 & 6 & 2.913109\\
    12 & 6 & 2.724988\\
    13 & 6 & 2.571162\\
    13 & 7 & 2.933424\\
    17 & 10 & 3.118090\\
\hline
\end{tabular}
\\
\vspace{5pt}
Values of exponents for~$K$
\end{minipage}
\hspace{5pt}
\begin{minipage}{0.4\textwidth}
\centering
\begin{tabular}{|c|c|c|}
\hline
$k$ & $s_k$ & $\Delta_k$ \\
\hline
4 & 3.5 & 0.849408 \\
4 & 3.55 & 0.819175\\
5 & 4 & 1.438656\\
6 & 5 & 1.724697\\
8 & 7 & 2.310599\\
10 & 8 & 3.414354 \\
\hline
\end{tabular}
\\
\vspace{5pt}
Values of exponents for~$L$
\end{minipage}
    \caption{Admissible exponents for \autoref{EstGlGk} used with the sets of \autoref{SetKL}}
    \label{TabAdmExpK7andK4}
\end{table}

\begin{table}[H]
    \centering
    \begin{minipage}{0.4\textwidth}
\centering
       \begin{tabular}{|c|c|c|}
\hline
$k$ & $s_k$ & $\Delta_k$ \\
\hline
4 & 3.3 & 0.977254 \\
4 & 3.35 & 0.944345\\
8 & 6 & 2.856656\\
8 & 7 & 2.310599\\
9 & 9 & 2.142503\\
10 & 8 & 3.414354\\
10 & 9 & 2.894571\\
11 & 9 & 3.693774\\
12 & 10 & 3.972927\\
13 & 11 & 4.252162\\
\hline
\end{tabular}
\\
\vspace{5pt}
Values of exponents for~$\tilde{K}$
\end{minipage}
\hspace{5pt}
\begin{minipage}{0.4\textwidth}
\centering
    \begin{tabular}{|c|c|c|}
\hline
$k$ & $s_k$ & $\Delta_k$ \\
\hline
5 & 5 & 0.92508 \\
5 & 6 & 0.541755\\
6 & 7 & 0.850572\\
7 & 7 & 1.541089\\
7 & 8 & 1.152632\\
17 & 40 & 0.317928\\
17 & 44 & 0.192585\\
\hline
\end{tabular}
\\
\vspace{5pt}
Values of exponents for~$\tilde{L}$
\end{minipage}
    \caption{Admissible exponents for \eqref{EqGKt} and \eqref{EqGLt}}
    \label{TabAdmExpGKtGLt}
\end{table}

\begin{table}[H]
    \centering
    \begin{minipage}{0.45\textwidth}
\centering
       \begin{tabular}{|c|c|c|}
\hline
$k$ & $s_k$ & $\Delta_k$ \\
\hline
4 & 1.65 & 2.35 \\
4 & 1.7 & 2.3\\
7 & 5 & 2.589117\\
8 & 5 & 2.589117\\
11 & 8 & 4.261426\\
12 & 8 & 5.132275\\
\hline
\end{tabular}
\\
\vspace{5pt}
with $s_4= \frac{1}{1-(1/5+1/5+1/8+1/8)}=\frac{5}{3}$
\end{minipage}
\hspace{1pt}
\begin{minipage}{0.5\textwidth}
\centering
    \begin{tabular}{|c|c|c|}
\hline
$k$ & $s_k$ & $\Delta_k$ \\
\hline
5 & 4 & 1.405522\\
5 & 5 & 0.925080\\
6 & 4 & 2.333333\\
9 & 6 & 3.738308\\
10 & 7 & 3.996716\\
13 & 9 & 5.385580\\
17 & 12 & 7.023320\\
\hline
\end{tabular}
\\
\vspace{5pt}
and $s_5= \frac{1}{1-(1/4+1/6+1/7+1/9+1/12)}=4.064516\dots$
\end{minipage}
    \caption{Admissible exponents for \eqref{EqGKGL}}
    \label{TabAdmExpGKGL}
\end{table}

\begin{table}[H]
    \centering
\begin{minipage}{0.23\textwidth}
\centering
\begin{tabular}{|c|c|c|}
\hline
$k$ & $s_k$ & $\Delta_k$\\
\hline
4 & 4.95 & 0.208007\\
4 & 5 & 0.187681 \\
6 & 8 &0.560413 \\
7 & 10 & 0.605568\\
9 & 13 & 0.895953\\
15 & 24 & 1.286149 \\
17 & 28 & 1.360721 \\
20 & 33 & 1.632692 \\
21 & 34 & 1.861960 \\
22 & 36 & 1.895913 \\
24 & 39 & 2.13053 \\
28 & 46 & 2.432087 \\
\hline
\end{tabular}
\end{minipage}
\hspace{15pt}
\begin{minipage}{0.23\textwidth}
\centering
\begin{tabular}{|c|c|c|}
\hline
$k$ & $s_k$ & $\Delta_k$\\
\hline
29 & 48 & 2.46681 \\
32 & 53 & 2.734333  \\
34 & 56 & 2.966921 \\
35 & 57 & 3.167478 \\
37 & 61 & 3.233989 \\
39 & 64 & 3.46723 \\
40 & 66 & 3.501023 \\
41 & 67 & 3.701116 \\
43 & 71 & 3.767860 \\
45 & 73 & 4.170490 \\
47 & 77 & 4.234431 \\
\hline
\end{tabular}
\end{minipage}
\hspace{15pt}
\begin{minipage}{0.23\textwidth}
\centering
\begin{tabular}{|c|c|c|}
\hline
$k$ & $s_k$ & $\Delta_k$\\
\hline
51 & 84 & 4.534192 \\
54 & 89 & 4.800804 \\
58 & 95 & 5.267750 \\
61 & 100 & 5.534095 \\
63 & 103 & 5.768085 \\
64 & 105 & 5.800461 \\
65 & 107 & 5.833356 \\
66 & 108 & 6.034346 \\
69 & 113 & 6.300573 \\
70 & 115 & 6.333120 \\
71 & 116 & 6.534702 \\
\hline
\end{tabular}
\end{minipage}
    \caption{Admissible exponents for \eqref{EqR31} used with the sets of \autoref{SetK4again}}
    \label{TabAdmExpR31K4}
\end{table}

\begin{table}[H]
    \centering
\begin{minipage}{0.23\textwidth}
\centering
\begin{tabular}{|c|c|c|}
\hline
$k$ & $s_k$ & $\Delta_k$\\
\hline
5 & 7 & 0.272729 \\
8 & 12 &  0.65612\\
10 & 15 & 0.930371 \\
11 & 17 & 0.967484 \\
12 & 18 &  1.189556\\
13 & 20 & 1.21907 \\
14 & 23 & 1.084921 \\
16 & 24 & 1.688314 \\
16 & 25 & 1.495533 \\
18 & 29 & 1.560877 \\
19 & 31 &  1.59611 \\
23 & 37 & 2.098122 \\
25 & 41 & 2.164099  \\
\hline
\end{tabular}
\end{minipage}
\hspace{15pt}
\begin{minipage}{0.23\textwidth}
\centering
\begin{tabular}{|c|c|c|}
\hline
$k$ & $s_k$ & $\Delta_k$\\
\hline
26 & 42 & 2.366005 \\
27 & 44 & 2.398294 \\
30 & 48 & 2.838828 \\
31 & 51 & 2.699544 \\
33 & 54 & 2.933149 \\
36 & 59 & 3.200317 \\
38 & 62 & 3.434398 \\
42 & 68 & 3.90423 \\
44 & 72 & 3.967842 \\
46 & 75 & 4.20209 \\
48 & 78 & 4.436808 \\
49 & 80 & 4.468516 \\
50 & 81 & 4.67194 \\
\hline
\end{tabular}
\end{minipage}
\hspace{15pt}
\begin{minipage}{0.23\textwidth}
\centering
\begin{tabular}{|c|c|c|}
\hline
$k$ & $s_k$ & $\Delta_k$\\
\hline
52 & 85 &  4.73497\\
53 & 86 & 4.937994 \\
55 & 89 & 5.173295 \\
56 & 91 & 5.173295 \\
57 & 93 & 5.235581 \\
59 & 96 & 5.470194 \\
60 & 98 & 5.501854 \\
62 & 101 & 5.736327 \\
67 & 109 & 6.237108 \\
68 & 110 & 6.441644 \\
72 & 117 & 6.737974 \\
73 & 119 & 6.769141 \\
\hline
\end{tabular}
\end{minipage}
    \caption{Admissible exponents for \eqref{EqR31} used with the sets of \autoref{SetK5}}
    \label{TabAdmExpR31K5}
\end{table}

\phantomsection
\addcontentsline{toc}{section}{References}

\bibliographystyle{alpha}
\bibliography{Refs}

@article{Bou16x,
  title = {On the {V}inogradov mean value},
  volume = {296},
  ISSN = {1531-8605},
  url = {http://dx.doi.org/10.1134/S0081543817010035},
  DOI = {10.1134/s0081543817010035},
  number = {1},
  journal = {Proc. Steklov Inst. Math. },
  publisher = {Pleiades Publishing Ltd},
  author = {Bourgain,  J.},
  year = {2017},
  month = Jan,
  pages = {30–40}
}

@article{Brudern1991,
author = {Br{\"{u}}dern,  J.},
journal = {Math. Scand.},
number = {1},
pages = {27-45},
title = {Ternary additive problems of {W}aring's type.},
url = {http://eudml.org/doc/167144},
volume = {68},
year = {1991},
}

@article{BruedernWooley2000,
  author  = {Br{\"u}dern, J. and Wooley, T. D.},
  title   = {On {W}aring's problem: two cubes and seven biquadrates},
  journal = {Tsukuba J. Math.},
  volume  = {24},
  number  = {2},
  year    = {2000},
  pages   = {387--417}
}

@article{Brudern2019,
  title = {Even ascending powers},
  author = {Br{\"{u}}dern,  J.},
  year = {2019},
  volume = {118},
  ISSN = {1730-6299},
  url = {http://dx.doi.org/10.4064/bc118-13},
  DOI = {10.4064/bc118-13},
  journal = {Banach Center Publ.},
  publisher = {Institute of Mathematics,  Polish Academy of Sciences},
  pages = {193–229}
}

@article{BrudernWooley2023Waringlarge,
url = {https://doi.org/10.1515/crelle-2023-0072},
title = {On {W}aring’s problem for larger powers},
author = {Br{\"{u}}dern,  J. and Wooley, T. D. },
pages = {115--142},
volume = {2023},
number = {805},
journal = {J. Reine Angew. Math.},
doi = {doi:10.1515/crelle-2023-0072},
year = {2023},
lastchecked = {2025-10-31}
}

@article{BrudernWooley2023Freiman,
   title={On {W}aring’s problem: Beyond {F}reĭman’s theorem},
   volume={109},
   ISSN={1469-7750},
   url={http://dx.doi.org/10.1112/jlms.12820},
   DOI={10.1112/jlms.12820},
   number={1},
   journal={Jour. London Math. Soc.},
   publisher={Wiley},
   author={Brüdern, J. and Wooley, T. D.},
   year={2023},
   month=oct }

@article{BrudernWooley2024SmoothMajor,
  title = {Estimates for Smooth {W}eyl Sums on Major Arcs},
  volume = {2024},
  ISSN = {1687-0247},
  url = {http://dx.doi.org/10.1093/imrn/rnae252},
  DOI = {10.1093/imrn/rnae252},
  number = {24},
  journal = {Int. Math. Res. Not. IMRN},
  publisher = {Oxford University Press (OUP)},
  author = {Br{\"{u}}dern,  J. and Wooley,  T. D},
  year = {2024},
  month = nov,
  pages = {14662–14688}
}

@article{Ford1995UnlikePowersI,
  author    = {K. B. Ford},
  title     = {The representation of numbers as sums of unlike powers},
  journal   = {J. London Math. Soc. (2)},
  year      = {1995},
  doi       = {10.1112/jlms/51.1.14},
  publisher = {London Mathematical Society}
}

@article{Ford1996UnlikePowersII,
  author    = {K. B. Ford},
  title     = {The representation of numbers as sums of unlike powers. {II}},
  journal   = {J. American Math. Soc.},
  year      = {1996},
  doi       = {10.1090/S0894-0347-96-00193-2},
  publisher = {American Mathematical Society}
}

@article{Freiman1949bVariantWaring,
  author    = {Fre{\u{\i}}man, G. A. },
  title     = {On one variant of {W}aring’s problem},
  journal   = {Vestnik Moskov. Univ. Ser. I Mat. Mekh.},
  year      = {1949},
  note      = {In Russian},
}

@article{Freiman1949aClaim,
  author  = {Fre{\u{\i}}man, G. A.},
  title   = {Solution of {W}aring's problem in a new form},
  journal = {Uspekhi Mat. Nauk},
  volume  = {4},
  number  = {1(29)},
  year    = {1949},
  pages   = {193}
}

@article{KLX2022,
  author    = {C. I. Kuan and D. Lesesvre and X. Xiao},
  title     = {Sums of even ascending powers},
  journal   = {arXiv},
  year      = {2020},
  note      = {Submitted for publication},
  url       = {https://arxiv.org/abs/2001.02429}
}

@phdthesis{Kufner2025,
  author       = {Küfner, T. S.},
  title        = {A real effective version of the Freĭman-Scourfield theorem},
  school       = {Georg-August-Universität Göttingen},
  year         = {2025},
  type         = {Ph.D. thesis},
  url          = {http://dx.doi.org/10.53846/goediss-11614},
}

@article{LiuZhao2021_RepresentationUnlikePowers,
  author    = {J. Liu and L. Zhao},
  title     = {Representation by sums of unlike powers},
  journal   = {J. Reine Angew. Math.},
  year      = {2021},
  doi       = {10.1515/crelle-2020-0001},
}

@article{Roth1951AdditiveNumberTheory,
  author    = {K. F. Roth},
  title     = {A problem in additive number theory},
  journal   = {Duke Math. J.},
  year      = {1951},
  doi       = {10.1215/S0012-7094-51-01868-8},
  publisher = {Duke University Press}
}

@article{Scourfield1963GeneralizationWaring,
  author    = {E. J. Scourfield},
  title     = {A generalization of {W}aring’s problem},
  journal   = {Proc. Edinburgh Math. Soc.},
  year      = {1963},
  doi       = {10.1017/S0013091500011324},
  publisher = {Cambridge University Press}
}

@article{Thanigasalam1983CertainAdditiveRepresentations,
  author    = {K. Thanigasalam},
  title     = {On certain additive representations of integers},
  journal   = {Portugaliae Math.},
  year      = {1983/84},
  publisher = {Sociedade Portuguesa de Matemática},
  url       = {http://eudml.org/doc/115560}
}

@book{vaughan1981hardy,
  author    = {R. C. Vaughan},
  title     = {The Hardy-Littlewood Method},
  publisher = {Cambridge University Press},
  year      = {1997},
  edition    = {Second},
  address   = {Cambridge}
}

@article{Vaughan1989,
  title = {A new iterative method in {W}aring’s problem},
  volume = {162},
  ISSN = {0001-5962},
  url = {http://dx.doi.org/10.1007/BF02392834},
  DOI = {10.1007/bf02392834},
  number = {0},
  journal = {Acta Math.},
  publisher = {International Press of Boston},
  author = {Vaughan,  R. C.},
  year = {1989},
  pages = {1–71}
}

@article{VaughanWooley1995,
  title = {Further improvements in {W}aring’s problem, {I}},
  volume = {174},
  ISSN = {0001-5962},
  url = {http://dx.doi.org/10.1007/BF02392467},
  DOI = {10.1007/bf02392467},
  number = {2},
  journal = {Acta Math.},
  publisher = {International Press of Boston},
  author = {Vaughan,  R. C. and Wooley,  T. D.},
  year = {1995},
  pages = {147–240}
}

@article{VaughanWooley1994,
  title = {Further improvements in {W}aring’s problem,  {II}: {S}ixth powers},
  volume = {76},
  ISSN = {0012-7094},
  url = {http://dx.doi.org/10.1215/S0012-7094-94-07626-6},
  DOI = {10.1215/s0012-7094-94-07626-6},
  number = {3},
  journal = {Duke Math. J.},
  publisher = {Duke University Press},
  author = {Vaughan,  R. C. and Wooley,  T. D.},
  year = {1994},
  month = dec 
}

@article{VaughanWooley1993,
  author  = {Vaughan, R. C. and Wooley, T. D.},
  title   = {Further improvements in {W}aring's problem, {III}: {E}ighth powers},
  journal = {Philos. Trans. Roy. Soc. London Ser. A},
  year    = {1993},
  pages   = {385--396}
}

@article{VaughanWooley2000,
author = {Vaughan, R. C. and Wooley, T. D.},
journal = {Acta Arith.},
language = {eng},
number = {3},
pages = {203-285},
title = {Further improvements in {W}aring's problem, {IV}: {H}igher powers},
url = {http://eudml.org/doc/207431},
volume = {94},
year = {2000}
}

@article{Wooley1991,
author = {Wooley,  T. D.},
journal = {J. Reine Angew. Math.},
pages = {141-198},
title = {On simultaneous additive equations {II.}},
url = {http://eudml.org/doc/153348},
volume = {419},
year = {1991},
}

@article{Wooley2015,
  title = {Sums of three cubes,  {II}},
  volume = {170},
  ISSN = {1730-6264},
  url = {http://dx.doi.org/10.4064/aa170-1-6},
  DOI = {10.4064/aa170-1-6},
  number = {1},
  journal = {Acta Arith.},
  publisher = {Institute of Mathematics,  Polish Academy of Sciences},
  author = {Wooley,  T. D.},
  year = {2015},
  pages = {73–100}
}

@article{Wooley2018,
  title = {Nested efficient congruencing and relatives of {V}inogradov’s mean value theorem},
  volume = {118},
  ISSN = {1460-244X},
  url = {http://dx.doi.org/10.1112/plms.12204},
  DOI = {10.1112/plms.12204},
  number = {4},
  journal = {Proc. London Math. Soc.},
  publisher = {Wiley},
  author = {Wooley,  T. D.},
  year = {2018},
  month = Oct,
  pages = {942–1016}
}

\adresseauteur{Lo\"{i}c Mermet}{Institut Camille Jordan CNRS UMR 5208}{Universit\'e Jean Monnet, CNRS, Centrale Lyon, INSA Lyon, Universit\'e Claude Bernard Lyon 1}{20 rue Annino, 42023 Saint-\'Etienne}{France}{loic.mermet@univ-st-etienne.fr}
\end{document}